\documentclass[a4paper,12pt]{article}

\usepackage{amsmath,amssymb,amsthm,mathtools}
\usepackage[T1]{fontenc}
\usepackage[utf8]{inputenc}
\usepackage{lmodern}
\usepackage[numbers,sort&compress]{natbib}
\usepackage[inline]{enumitem}
\usepackage[title]{appendix}
\usepackage{hyperref}
\hypersetup{plainpages=false,pdfstartview={FitV null},colorlinks=true,
  linkcolor=blue,citecolor=blue,urlcolor=blue,
  pdftitle={An Elementary Approach to Depoissonization},
  pdfauthor={Vytas Zacharovas}}

\DeclareRobustCommand{\Stirling}{\genfrac\{\}{0pt}{}}
\newtheorem{thm}{Theorem}[section]
\newtheorem{cor}[thm]{Corollary}
\newtheorem{lem}[thm]{Lemma}
\newtheorem{prop}[thm]{Proposition}
\theoremstyle{remark}
\newtheorem{rem}[thm]{Remark}
\newtheorem{exm}[thm]{Example}

\let\originalappendices\appendices
\renewcommand{\appendices}{\originalappendices
  \setcounter{equation}{0}
  \renewcommand{\theequation}{\Alph{section}\arabic{equation}}}
\title{An Elementary Approach to Depoissonization}
\author{Vytas Zacharovas\thanks{Institute of Computer Science,
Vilnius University, Naugarduko g.~24, Vilnius, 03225, Lithuania.
E-mail: \mbox{\href{mailto:vytas.zacharovas@mif.vu.lt}{vytas.zacharovas@mif.vu.lt}}}\thanks{Institute of Statistical Science,
Academia Sinica, No.\ 128, Academia Road, Section~2, Nankang,
Taipei, 115201, Taiwan.}}
\date{}

\begin{document}
\maketitle
\begin{abstract}
We study depoissonization for sequences with entire exponential
generating functions. We establish explicit finite-order remainder
bounds for Poisson--Charlier approximations, controlled by
Poisson-weighted averages of the absolute values of higher-order
forward differences of the coefficient sequence. This complements
classical analytic depoissonization, which instead relies on
complex-plane growth estimates: classical methods are natural when
such estimates are available, whereas the present approach is useful
when the relevant finite differences can be bounded, without requiring
estimates at nonreal arguments.

We also study the inverse problem of deriving asymptotic expansions of
the Poisson transform at large positive arguments from the coefficient
sequence. This leads to finite-difference analogues of Ramanujan's
expansion. Applications to combinatorics and probability are presented.
\end{abstract}

\noindent\emph{MSC 2020 Subject Classifications}: Primary 41A60;
secondary 05A15, 60C05.

\noindent\emph{Key words}: Depoissonization,  Ramanujan's expansion, Charlier polynomials,  Operator bounds, Tauberian theory.

\section{Introduction}
\subsection{Background and motivation}
Let $(A_m)_{m\geqslant 0}$ be a sequence of real numbers such that the
exponential generating series
\[
\sum_{m=0}^{\infty} A_m\frac{x^m}{m!}
\]
converges for every $x\in\mathbb{R}$. Define the corresponding Poisson
transform by
\begin{equation}
\label{eq:poisson-transform}
f(x):=e^{-x}\sum_{m=0}^{\infty} A_m\frac{x^m}{m!}.
\end{equation}
Since the power series on the right-hand side has infinite radius of
convergence, the function $f$ is infinitely differentiable on
$\mathbb{R}$, and the series may be differentiated term by term any
finite number of times. Equivalently, the same power series defines an entire function of a
complex variable, although the estimates developed below use only
real arguments. We consider situations in which either a closed-form
expression for \(f(x)\) is available or \(f(x)\) can be estimated for
positive \(x\), while the asymptotic behavior of \(A_m\) is unknown. For
\(x>0\), the numbers
\[
e^{-x}\frac{x^m}{m!},\qquad m=0,1,2,\ldots,
\]
form the probability mass function of a Poisson random variable with
parameter \(x\). In particular, if \(\eta(n)\) is a Poisson random variable with
parameter \(n\), then
\[
f(n)=\mathbb{E}A_{\eta(n)}.
\]
Since \(\eta(n)\) has mean and variance \(n\), for every
\(T_n\to\infty\) the probability that
\[
|\eta(n)-n|>T_n\sqrt n
\]
tends to zero. Hence, for a bounded sequence \((A_m)\),
\[
f(n)-A_n
=
\sum_{\substack{m\geqslant0\\ |m-n|\leqslant T_n\sqrt n}}
e^{-n}\frac{n^m}{m!}\bigl(A_m-A_n\bigr)
+o(1).
\]
Consequently, \(f(n)=A_n+o(1)\) if the sequence varies sufficiently
slowly on the natural \(\sqrt n\)-scale. For example, under the same
boundedness assumption, the condition
\[
|A_{m+1}-A_m|=o\!\left(m^{-1/2}\right),
\qquad m\to\infty,
\]
is sufficient: by telescoping, one can choose \(T_n\to\infty\)
sufficiently slowly so that \(A_m-A_n=o(1)\) uniformly for
\(|m-n|\leqslant T_n\sqrt n\). Theorem~\ref{thm:main-depoissonization}
and the discussion following it give a quantitative version of this
principle without the boundedness assumption. Deriving asymptotic estimates for \(A_n\) from the Poisson-smoothed
values \(f(n)\), under suitable regularity assumptions on the sequence
\((A_m)\), is a typical Tauberian problem; see \citet{korevaar_2004}.

The  approximation \(A_n\approx f(n)\) can be also  regarded as the first-term asymptotic expansion 
of \(A_n\) in terms of \(f(x)\). 
To obtain a more accurate asymptotic expression, we expand \(f(x)\) 
into its Taylor series at the point \(x = n\):
\[
\sum_{\ell=0}^\infty \frac{A_\ell}{\ell!}\,x^\ell
= e^x f(x)
= e^x \sum_{m=0}^\infty \frac{f^{(m)}(n)}{m!}\,(x - n)^m.
\]
Comparing the coefficients of \(x^n\) on both sides yields
\begin{equation}
\label{eq:charlier-poisson-expansion}
A_n = \sum_{m=0}^\infty \frac{f^{(m)}(n)}{m!}\,\tau_m(n),
\end{equation}
where
\begin{equation}
\label{eq:tau-coefficient-definition}
\tau_m(n) = n!\,[x^n]\bigl(e^x (x - n)^m\bigr),
\end{equation}
and the notation \([x^n]F(x)\) denotes the coefficient of \(x^n\) 
in the power-series expansion of the function \(F(x)\) around the origin:
\[
F(x) = \sum_{k=0}^\infty a_k x^k
\quad \Longrightarrow \quad
[x^n]F(x) = a_n.
\]
We refer to~\eqref{eq:charlier-poisson-expansion} as the
\emph{Poisson--Charlier expansion}, following the terminology
introduced in~\citet{hwang_fuchs_zacharovas_2010}.
One can easily see that if \(f\) is an entire function, the   identity \eqref{eq:charlier-poisson-expansion}   holds for all \(n\).
The functions \(\tau_s(x)\) are polynomials in \(x\) of degree at most
\(\lfloor s/2\rfloor\) and may alternatively be defined by the
generating function
\begin{equation}
	\label{eq:tau-generating-function-intro}
	\sum ^{\infty }_{m=0}\frac{\tau_{m}\left( x\right) }{m!}z^m=\bigl((1+z)e^{-z}\bigr)^x,
	\qquad |z|<1.
\end{equation}
The polynomials \(\tau_m\) are the diagonal specializations of the
Charlier polynomials introduced later in~\eqref{eq:charlier-coefficient-definition}; more
precisely,
\[
\tau_m(n)=n^m C_m(n,n),
\qquad m\geqslant0,\quad n\geqslant1.
\]
The first several of these polynomials are as follows:
\[
\renewcommand{\arraystretch}{1.1}
\begin{array}{|c|c|c|c|c|c|c|}
\hline
\tau_{0}(x) & \tau_{1}(x) & \tau_{2}(x) & \tau_{3}(x) & 
\tau_{4}(x) & \tau_{5}(x) & \tau_{6}(x) \\
\hline
1 & 0 & -x & 2x & 
3x^{2}-6x & 24x-20x^{2} & -15x^{3}+130x^{2}-120x \\
\hline
\end{array}
\]
More information on the properties of these polynomials is provided in the Appendix~ \ref{app:polynomial-properties}.  It is known that, under additional conditions on \(f\), the
expansion~\eqref{eq:charlier-poisson-expansion} is also an asymptotic
expansion: as \(n\to\infty\), finite truncations containing
successively more terms provide approximations of \(A_n\) with
correspondingly smaller asymptotic remainders.
The following theorem is  based on results from \citet{jacquet_szpankowski_1998}, with two key differences. First, we impose less general conditions on the generating function to convey the core idea of the result more transparently, avoiding technical details that might obscure the main argument. Second, the asymptotic term is expressed using the polynomials $\tau_m$, following the approach of \citet{hwang_fuchs_zacharovas_2010, fuchs_hwang_zacharovas_2014}. This representation is chosen for its compatibility with the subsequent developments in this paper.
\begin{thm}[\citet{jacquet_szpankowski_1998}]
\label{thm:classical-depoissonization}
Let \(f\) denote the entire extension to \(\mathbb{C}\) of the Poisson
transform defined in~\eqref{eq:poisson-transform}, so that
\[
f(z)=e^{-z}\sum_{m=0}^{\infty}\frac{A_m}{m!}z^m,
\qquad z\in\mathbb{C}.
\]
Assume that, for some fixed \(\varepsilon>0\):
\begin{enumerate*}[label=(\roman*)]
\item in the sector \( |\arg z| \leqslant \varepsilon \), one has \( f(z) = O(|z|^{\alpha}) \) as \( |z| \to \infty \);
\item in the complementary region \( \varepsilon \leqslant |\arg z| \leqslant \pi \), the function \( e^{z} f(z) \) satisfies 
\( e^{z} f(z) = O(e^{(1 - \delta)|z|}) \) for some \( \delta > 0 \).
\end{enumerate*}
Then, for every integer \( k \geqslant 1 \),
\[
A_n = \sum_{0 \leqslant j < 2k} \frac{f^{(j)}(n)}{j!}\,\tau_j(n)
      + O\!\left(n^{\alpha - k}\right),
\]
where \( \tau_j(n) \) are defined in~\eqref{eq:tau-coefficient-definition}.
\end{thm}

Theorem~\ref{thm:classical-depoissonization} follows the classical depoissonization approach: it
derives an asymptotic expansion of the coefficients \(A_n\) from estimates
for \(f(z)\) throughout suitable regions of the complex plane. One of the
principal aims of this paper is to develop an alternative approach that
avoids such complex-analytic information. Instead, our results rely on
upper bounds for finite differences of the coefficient sequence, together
with information about \(f(x)\) for positive real values of the argument.
This approach is not intended to replace the classical method or to be
regarded as inherently stronger or weaker; rather, the two methods are
complementary. The present approach is designed primarily for situations
in which the classical complex-analytic conditions are difficult to verify,
for example when the analytic continuation of \(f\) has a complicated
form, is defined only implicitly, or is otherwise difficult to estimate
for complex values of the argument.

Throughout the paper, unless stated otherwise, $(A_m)_{m\geqslant 0}$
denotes a sequence of real numbers whose exponential generating series
converges for every $x\in\mathbb{R}$, and $f$ denotes the corresponding
Poisson transform defined by~\eqref{eq:poisson-transform}. The series has infinite
radius of convergence and may therefore be differentiated term by term
any finite number of times on $\mathbb{R}$.

We use $\Delta$ to denote the forward difference operator.
For a function $h$ defined on $\mathbb{N}_0$ or on $[0,\infty)$, we set
\[
(\Delta h)(x):=h(x+1)-h(x),
\]
and write $\Delta^k$ for its $k$-fold iteration, with
$\Delta^0 h=h$. In particular, when the sequence $(A_m)$ is regarded
as a function on $\mathbb{N}_0$,
\[
\Delta A_m=A_{m+1}-A_m
\]
and, for example,
\[
\Delta^2A_m
=A_{m+2}-2A_{m+1}+A_m.
\]

Repeated differentiation of~\eqref{eq:poisson-transform} then gives, for every
integer $k\geqslant 0$,
\begin{equation}
\label{eq:derivative-difference-correspondence}
f^{(k)}(x)
=
e^{-x}\sum_{m=0}^{\infty}
\frac{\Delta^k A_m}{m!}x^m.
\end{equation}
Thus, differentiation of the Poisson transform corresponds to forward
differencing of the coefficient sequence of $e^x f(x)$. This
correspondence will be used repeatedly throughout the paper. For any entire function \(f\) of the form~\eqref{eq:poisson-transform} and any
\(r\geqslant0\), define
\begin{equation}
\label{eq:E-definition}
E(f;r):=e^{-r}\sum_{m=0}^{\infty}\frac{|A_m|}{m!}r^m.
\end{equation}
We use the same notation when the Poisson coefficients \(A_m\) are complex,
with \(|A_m|\) denoting the complex modulus.
If a function \(f\) depends on several variables, we indicate the variable with respect to which the operator \(E\) is taken by writing it as a subscript of \(E\). For example, \(E_x(f(ax);r)\) denotes \(E(h;r)\), where \(h(x)=f(ax)\).
 An example of our results is the following theorem.
 
\begin{thm}
\label{thm:main-depoissonization}
Let $(A_m)_{m\geqslant 0}$ be a sequence of real numbers such that
\[
\sum_{m=0}^{\infty} A_m\frac{x^m}{m!}
\]
converges for every $x\in\mathbb{R}$, and define
\[
f(x)=e^{-x}\sum_{m=0}^{\infty} A_m\frac{x^m}{m!}.
\]
Then, for all integers $n\geqslant 2$ and $N\geqslant 0$, one has
\[
\left|
A_n-\sum_{m=0}^{N}\frac{f^{(m)}(n)}{m!}\tau_m(n)
\right|
\leqslant
33\,n^{\frac{N+1}{2}}
E\!\left(f^{(N+1)};n\right).
\]
\end{thm}

By~\eqref{eq:E-definition}, for \(r\geqslant0\),
\[
E\!\left(f^{(N+1)}; r\right) 
   = e^{-r} \sum_{m=0}^\infty \frac{|\Delta^{N+1} A_m|}{m!}\, r^m,
\]
so the theorem shows that the coefficient \(A_n\) can be approximated by a
partial sum of its Poisson--Charlier expansion~\eqref{eq:charlier-poisson-expansion},
provided that a suitable Poisson-weighted average of the finite differences 
\(|\Delta^{N+1} A_m|\) near \(m = n\) is 
\(o\!\left(n^{-(N+1)/2}\right)\). For example, if we take \(N=0\) and assume that 
\(\Delta A_m = o\!\left(m^{-1/2}\right)\) as \(m\to\infty\), then a routine estimate yields 
\(E(f';R) = o\!\left(R^{-1/2}\right)\). In particular, the right-hand side of 
the inequality in Theorem~\ref{thm:main-depoissonization} tends to zero, which gives
\(A_n - f(n) = o(1)\). A classical result of Hardy and Littlewood 
(see Theorem~156 in \citet{hardy_1992} or Theorem~9.1 in
\citet{korevaar_2004}) further shows that if, in addition, \(f(n) \to C\) as 
\(n \to \infty\), then even the weaker condition
\(
|\Delta A_m| = O\!\left(m^{-1/2}\right)
\)
in place of \(o\!\left(m^{-1/2}\right)\) already suffices to guarantee that 
\(A_n \to C\).

A particularly simple situation arises when the finite differences of order
\(N+1\) all have the same sign; that is, when either
\(\Delta^{N+1} A_m \geqslant 0\) for all \(m \geqslant 0\) or
\(\Delta^{N+1} A_m \leqslant 0\) for all \(m \geqslant 0\).
Since the weights \(e^{-r}r^m/m!\) are nonnegative for \(r\geqslant0\),
the fixed-sign assumption and the derivative--difference identity give
\[
\begin{aligned}
E\!\left(f^{(N+1)};r\right)
&=
\left|
e^{-r}\sum_{m=0}^{\infty}
\frac{\Delta^{N+1}A_m}{m!}r^m
\right| \\
&=
\bigl|f^{(N+1)}(r)\bigr|,
\qquad r\geqslant0.
\end{aligned}
\]
At \(r=0\), both sides are equal to
\(\lvert\Delta^{N+1}A_0\rvert\), so strict positivity of \(r\) is not
required.
Consequently, Theorem~\ref{thm:main-depoissonization} gives an approximation with an \(o(1)\) remainder
whenever 
\(
\bigl| f^{(N+1)}(n) \bigr| = o\!\left(n^{-(N+1)/2}\right)
\)
as \(n \to \infty\).

In the simplest case \(N=0\), the requirement that all first-order finite
differences of \((A_m)\) have the same sign simply means that the sequence is
monotone. In this setting one obtains a sharper estimate, formulated as the
following corollary.

\begin{cor}
\label{cor:monotone-depoissonization}
If \((A_m)\) is monotone (either nondecreasing or nonincreasing), then the
inequality
\[
\bigl| A_n - f(R) \bigr|
   \;\leqslant\; \frac{e^R}{R^n}\, n! \,\bigl| f'(R) \bigr|
\]
holds for all \(R>0\) and \(n \geqslant 1\).
In particular, for \(R=n\) one has the sharper bound
\begin{equation}
\label{eq:monotone-diagonal-bound}
\bigl| A_n - f(n) \bigr| 
   \;\leqslant\; 2\sqrt{n}\,\bigl| f'(n) \bigr|.
\end{equation}
\end{cor}

\begin{rem}
The corollary shows that, for a monotone sequence, the
depoissonization error can be controlled solely by the first
derivative of the Poisson transform:
\[
|A_n-f(n)|
\leqslant
2\sqrt n\,|f'(n)|.
\]
Thus, in this setting, depoissonization reduces to estimating
\(f'(n)\).

A related monotonicity-and-sandwiching method was introduced in
Lemma~2.5 of \citet{johansson_1998} and generalized to generalized
Poisson generating functions by
\citet[Lemma~5.1]{bornemann_2025}. The following
standard-Poisson formulation captures their common underlying idea. Suppose that
\((A_n)\) is nonnegative, nondecreasing, and bounded above by a
constant \(C\). If \(\eta(\lambda)\) is Poisson with parameter
\(\lambda\), then
\[
\frac{f(\lambda)
      -C\mathbb P\!\left(\eta(\lambda)>n\right)}
     {\mathbb P\!\left(\eta(\lambda)\leqslant n\right)}
\leqslant A_n
\leqslant
\frac{f(\lambda)}
     {\mathbb P\!\left(\eta(\lambda)\geqslant n\right)}.
\]

Since \(A_n\) is independent of \(\lambda\), different parameters may
be used in the two bounds. Let
\[
\lambda_\pm=n\pm T_n\sqrt n,
\qquad
T_n\to\infty,
\qquad
T_n=o(\sqrt n).
\]
The relevant Poisson probabilities in the denominators then tend to
\(1\), while
\(\mathbb P(\eta(\lambda_-)>n)\to0\). Since \(0\leqslant f\leqslant C\),
the preceding inequalities reduce to
\[
f(\lambda_-)+o(1)
\leqslant A_n
\leqslant
f(\lambda_+)+o(1).
\]
Moreover, the mean-value theorem gives
\[
|f(\lambda_\pm)-f(n)|
\leqslant
T_n\sqrt n\,
\sup_{|x-n|\leqslant T_n\sqrt n}|f'(x)|.
\]
Consequently, a sufficient condition for Johansson's argument to yield
\(A_n=f(n)+o(1)\) is
\[
T_n\sqrt n\,
\sup_{|x-n|\leqslant T_n\sqrt n}|f'(x)|
\longrightarrow0.
\]

The factor \(\sqrt n\) reflects the Poisson standard-deviation scale
that also appears in Corollary~\ref{cor:monotone-depoissonization}. The distinction
is that Johansson's argument requires control at shifted arguments, or
uniformly on a growing neighborhood of \(n\), whereas
Corollary~\ref{cor:monotone-depoissonization} gives a direct estimate in terms of
\(f'(n)\) alone. The corollary also applies to both nondecreasing and
nonincreasing sequences and does not require nonnegativity or
boundedness.

\end{rem}

The monotone estimate also has a probabilistic application: it controls the
difference between fixed-index and poissonized expectations in terms of the
variance of the poissonized variable. This application is formulated in
Proposition~\ref{prop:poissonized-expectation-bound} in Section~\ref{sec:applications}.

The same monotonicity principle also yields a uniform approximation for
the Stirling numbers of the second kind; see
Proposition~\ref{prop:stirling-first-order} in
Section~\ref{sec:applications}.
 \subsection{Poissonization and Ramanujan's asymptotic expansion}
 A problem that may be regarded as the inverse of the depoissonization problem was studied by Ramanujan (see p.~58 of \citet{berndt_1985}) and later generalized by \citet{yu_2009}. Specifically, Ramanujan considered the task of determining the asymptotic
behavior of
\[
f(x)=e^{-x}\sum_{n=0}^{\infty}\frac{A_n}{n!}x^n,
\qquad x\to+\infty,
\]
assuming that the sequence \((A_n)\) is represented at the nonnegative
integers by a sufficiently regular function \(\varphi\), so that
\[
\varphi(n)=A_n,\qquad n=0,1,2,\ldots,
\]
and that the derivatives of \(\varphi\) satisfy suitable quantitative
growth conditions. The following formulation is Berndt's version of Ramanujan's
Entry~10, supplemented by the lower-bound hypothesis on \(G\)
appearing in Theorem~1 of \citet{yu_2009}.
\begin{thm}[A strengthened form of Ramanujan's Entry~10]
\label{thm:ramanujan}
Let $\varphi(x)$ denote a function of at most polynomial growth as $x\to+\infty$ (for real $x$).
Suppose there exist constants \(x_0>0\) and \(C\geqslant1\), and
a function \(G(x)\) of at most polynomial growth as
\(x\to+\infty\), such that, for every nonnegative integer \(m\)
and every \(x\geqslant x_0\), the derivative
\(\varphi^{(m)}(x)\) exists and satisfies
\begin{equation}\label{eq:ramanujan-derivative-bound}
  \left|\frac{\varphi^{(m)}(x)}{m!}\right|
  \leqslant G(x)\!\left(\frac{C}{x}\right)^{m}.
\end{equation}
Assume also that there exist constants \(c_0,c_1>0\) such that
\[
G(x)\geqslant c_0e^{-c_1\sqrt{x}}
\]
for all sufficiently large \(x\).
Then, for any fixed positive integer $M$,
\begin{equation}\label{eq:ramanujan-expansion}
\begin{aligned}
e^{-R}\sum_{m=0}^{\infty} \frac{\,\varphi(m)}{m!}R^{m}
  &= \varphi(R) \\
  &\quad+ \sum_{k=2}^{M}\ \sum_{n=k}^{2k-2}
      b_{k,n}\, R^{\,n-k+1}\, \frac{\varphi^{(n)}(R)}{n!} \\
  &\quad+ O\!\big(G(R)\,R^{-M}\big)
  \qquad (R\to+\infty).
\end{aligned}
\end{equation}
Here the coefficients $b_{k,n}$ are the nonnegative integers determined by
\begin{equation}\label{eq:ramanujan-coefficient-recurrence}
\begin{cases}
b_{k ,k}=1,\quad b_{k, n}=0, & \text{for } n<k \text{ or } n>2k-2,\\[3pt]
b_{k+1,\,n+1}=n\,b_{k,\,n-1}+(n-k+1)\,b_{k,\,n}, & \text{for } k\leqslant n\leqslant 2k-1 .
\end{cases}
\end{equation}
\end{thm}
The lower-bound hypothesis on \(G\) is absent from the formulation
of Entry~10 in \citet[pp.~57--58]{berndt_1985}, but is included in
Theorem~1 of \citet{yu_2009}. It allows a separate error term that
is exponentially small in \(\sqrt{x}\), arising in Berndt's proof,
to be absorbed into the stated \(O(G(x)x^{-M})\) remainder.

\citet{yu_2009} interpreted the left-hand side of the asymptotic
expansion~\eqref{eq:ramanujan-expansion} as the expectation
\(\mathbb{E}\varphi(\eta(R))\), where \(\eta(R)\) is a Poisson random
variable with parameter \(R\). Using this probabilistic interpretation,
he extended Ramanujan's result to expectations
\(\mathbb{E}\varphi(Y)\) for broader classes of convolution-type
families of random variables.

In the Poisson case, Yu's result also reduces the smoothness
requirements of Ramanujan's theorem. For a fixed truncation order, the
existence of all derivatives of \(\varphi\) can be replaced by the
existence of only finitely many derivatives
\(\varphi^{(m)}\), \(1\leqslant m\leqslant 2M\), together with the
corresponding bounds of the form~\eqref{eq:ramanujan-derivative-bound}. However, Yu's formulation still requires both \(\varphi\) and
\(G\) to have at most polynomial growth, together with the local
comparability condition on \(G\) stated in Theorem~2.

Theorem~\ref{thm:ramanujan-real-finite-difference} and Corollary~\ref{cor:ramanujan-derivative-form}
address a different limitation. They replace the collection of lower
derivative bounds by a condition only on the highest derivative involved
in the chosen truncation, while allowing substantially faster growth of
the interpolation function.

The applicability of Theorem~\ref{thm:ramanujan} to a given sequence \((A_n)\) therefore
depends on the construction of an interpolation \(\varphi\) for which
the required regularity and derivative-growth estimates can be
established. In applications, constructing and analysing such an
interpolation may be considerably more difficult than working with the
coefficient sequence itself. This motivates us to seek, as a first
step, a Ramanujan-type expansion at integer arguments that is formulated
directly in terms of the discrete sequence and does not require such an
interpolation.

The Poisson--Charlier expansion~\eqref{eq:charlier-poisson-expansion}
expresses \(A_n\) in terms of the derivatives of \(f\) at \(n\).
Since differentiation of the Poisson transform corresponds to forward
differencing of the coefficient sequence, it is natural to ask whether
this expansion can be formally inverted to express \(f(n)\) in terms of
the finite differences \(\Delta^m A_n\).

To see what such an inversion should look like, let
\(D_x=d/dx\). Then~\eqref{eq:charlier-poisson-expansion} may be written as
\[
A_n=
\left(
\sum_{m=0}^\infty \frac{\tau_m(n)}{m!}D_x^m
\right)
f(x)\Biggl|_{x=n}.
\]
By the derivative--difference correspondence
\eqref{eq:derivative-difference-correspondence},
\(f^{(j)}\) is the Poisson transform of the sequence
\((\Delta^jA_m)_{m\geqslant0}\). Applying
\eqref{eq:charlier-poisson-expansion} to this sequence therefore gives
\[
\Delta^jA_n=
\left(
\sum_{m=0}^\infty \frac{\tau_m(n)}{m!}D_x^m
\right)
D_x^j f(x)\Biggl|_{x=n}.
\]
Multiplying by \(\tau_j(-n)/j!\) and summing over \(j\), we obtain
\begin{equation}
\label{eq:formal-inverse-operator}
\sum_{j=0}^\infty \frac{\tau_j(-n)}{j!}\Delta^jA_n
=
\left(
\sum_{m=0}^\infty \frac{\tau_m(n)}{m!}D_x^m
\right)
\left(
\sum_{j=0}^\infty \frac{\tau_j(-n)}{j!}D_x^j
\right)
f(x)\Biggl|_{x=n}.
\end{equation}
By the generating function~\eqref{eq:tau-generating-function-intro},
\[
\left(
\sum_{m=0}^\infty\frac{\tau_m(n)}{m!}z^m
\right)
\left(
\sum_{j=0}^\infty\frac{\tau_j(-n)}{j!}z^j
\right)
=
\bigl((1+z)e^{-z}\bigr)^n
\bigl((1+z)e^{-z}\bigr)^{-n}
=1.
\]
Formally replacing \(z\) by \(D_x\) therefore makes the product of the
two differential-operator series in
\eqref{eq:formal-inverse-operator} the identity operator. This suggests
the formal inversion
\begin{equation}
\label{eq:ramanujan-inverted-series}
f(n)
=
\sum_{m=0}^{\infty}
\frac{\Delta^m A_n}{m!}\tau_m(-n).
\end{equation}
The passage from the reciprocal generating-function identity to the
differential-operator calculation above is purely formal: convergence
and the interchange or rearrangement of the resulting infinite series
require separate justification.

More fundamentally, the formal identity
\eqref{eq:ramanujan-inverted-series} cannot hold for an arbitrary
sequence. Its right-hand side depends only on
\(A_n,A_{n+1},\ldots\), whereas \(f(n)\) depends on the whole sequence.
For example, the identity generally fails when the sequence has finite
support. Proposition~\ref{prop:laplace-sequence-inversion} gives a
nontrivial class of sequences for which
\eqref{eq:ramanujan-inverted-series} does hold exactly.

The following theorem addresses a different question. It does not
assert convergence of the infinite series in~\eqref{eq:ramanujan-inverted-series} for
fixed \(n\). Instead, it gives a rigorous error bound for every finite
truncation, with the truncation order \(N\) fixed as \(n\to\infty\).
Thus, when the corresponding finite-difference remainder is sufficiently
small, the formal series may be interpreted as an asymptotic expansion.
\begin{thm}
\label{thm:ramanujan-integer-finite-difference}
Let $(A_m)_{m\geqslant 0}$ be a sequence of real numbers such that
\[
\sum_{m=0}^{\infty} A_m\frac{x^m}{m!}
\]
converges for every $x\in\mathbb{R}$, and define
\[
f(x)=e^{-x}\sum_{m=0}^{\infty} A_m\frac{x^m}{m!}.
\]
Then, for each fixed integer $N\geqslant 0$, uniformly for all integers
$n\geqslant 2$,
\[
\left|
f(n)-\sum_{m=0}^{N}
\frac{\Delta^m A_n}{m!}\tau_m(-n)
\right|
\leqslant
33(N+1)n^{\frac{N+1}{2}}
E\!\left(f^{(N+1)};n\right).
\]
\end{thm}
Theorem~\ref{thm:ramanujan-integer-finite-difference} gives an expansion of a particularly
simple form, but only at integer arguments. For arbitrary $R>0$,
Theorem~\ref{thm:reverse-depoissonization-arbitrary-R} gives an extension in which the
finite differences are taken at $n=\lfloor R\rfloor$ and the
coefficients involve Charlier polynomials depending separately on
$R$ and $n$. It is therefore not immediate that the simpler form of
Theorem~\ref{thm:ramanujan-integer-finite-difference} should survive the passage from $n$
to a real argument $R$.

Theorem~\ref{thm:ramanujan-real-finite-difference} shows that it does, provided that $(A_n)$ is the
restriction of a sufficiently regular function $\varphi$ satisfying
the derivative and growth conditions stated below. The expansion is
then expressed directly through $\Delta^s\varphi(R)$ and
$\tau_s(-R)$. These quantitative conditions, rather than the mere
existence of an interpolation, are essential.

\begin{thm} 
\label{thm:ramanujan-real-finite-difference}
Fix an integer $N\geqslant 0$. Suppose that the sequence $(A_n)_{n\geqslant 0}$
is the restriction to the nonnegative integers of a function
$\varphi\in C^{N+1}([0,\infty))$, so that
\[
A_n=\varphi(n),\qquad n=0,1,2,\ldots,
\]
and assume that its $(N+1)$-st derivative satisfies\begin{equation}
\label{eq:highest-derivative-bound}
\left|{\varphi^{(N+1)}(x)}\right|\leqslant \frac{G(x)}{x^{N+1}},	
\end{equation}
for all \(x\geqslant T\), where \(T>0\) is fixed and \(G(x)>0\) is a function such that
\begin{equation}
\label{eq:variation-condition}
{G(y)}\leqslant B\, {G(x)}e^{C\frac{|y-x|}{\sqrt{x}}} ,
\end{equation}
for all \(x,y\geqslant T\), with some constants \(B,C>0\). Then, as \(R\to+\infty\),
\[
e^{-R}\sum_{m=0}^\infty \frac{\varphi(m)}{m!}R^m
=\sum_{s=0}^{N}\frac{\Delta^s \varphi(R)}{s!}\,\tau_{s}(-R)
+O\!\left(R^{-\frac{N+1}{2}}G(R)\right)
+O\!\left(R^{T+(N+1)/2}e^{-R}\right).
\]
Here $\Delta$ is applied to the function $\varphi$ on $[0,\infty)$. 
The constants implied by the \(O(\cdot)\)-terms may depend on \(N\), \(B\), \(C\), \(T\), and on the values of \(\varphi\) on a fixed bounded interval, but are independent of \(R\).
\end{thm}
The condition \eqref{eq:variation-condition} can be compared to the corresponding condition of Yu (see Theorem 2 of \citet{yu_2009}) which required
\[
G(y)\leqslant  BG(x)
\]
for some fixed constant \(B\) and all \(|y-x|\leqslant \eta |x|\) with some fixed \(\eta\in(0,1)\). Our condition allows \(G\) to vary much faster  when \(|x-y|\gg \sqrt{x}\), while on the other hand Yu's result applies to a very wide range of other distributions, not only the Poisson distribution.

Theorem~\ref{thm:ramanujan-real-finite-difference} shows that a similar
expansion still holds when one retains only the upper bound for the
highest derivative and drops the corresponding bounds on the lower
derivatives.
\begin{cor}
\label{cor:ramanujan-derivative-form}
Under the assumptions of Theorem~\ref{thm:ramanujan-real-finite-difference}, we have
\[
e^{-R}\sum_{m=0}^\infty \frac{\varphi(m)}{m!}R^m
=\sum_{j=0}^{N}\frac{\varphi^{(j)}(R)}{j!}\,\rho_j(-R)
+O\!\left(R^{-\frac{N+1}{2}}G(R)\right)
+O\!\left(R^{T+(N+1)/2}e^{-R}\right),
\]
as \(R\to+\infty\),
where \(\rho_j(x)\) denote Mahler polynomials (see \citet[p.~254]{erdelyi_et_al_1981_vol3}), defined by their generating function
\begin{equation}
\label{eq:mahler-generating-function}
e^{-x(e^z-1-z)}=\sum_{j=0}^\infty \frac{\rho_j(x)}{j!}z^j.
\end{equation}
In particular, the polynomials \(\rho_j(-x)\) are of degree at most \(\lfloor j/2\rfloor\) in \(x\), with nonnegative coefficients.
The first eight of them are given below:
\[
\renewcommand{\arraystretch}{1.1}
\begin{array}{|c|c|c|c|c|c|c|c|}
\hline
\rho_{0}(-x) & \rho_{1}(-x) & \rho_{2}(-x) & \rho_{3}(-x) & 
\rho_{4}(-x) & \rho_{5}(-x) & \rho_{6}(-x) & \rho_{7}(-x) \\
\hline
1 & 0 & x & x &
3x^{2}+x & 10x^{2}+x & 15x^{3}+25x^{2}+x & 105x^{3}+56x^{2}+x \\
\hline
\end{array}
\]
The constants implied by both \(O(\cdot)\)-terms have the same
dependencies as in Theorem~\ref{thm:ramanujan-real-finite-difference}: they may depend on the
fixed parameters \(N,B,C,T\) and on the bounded-interval data for
\(\varphi\), but not on \(R\).
\end{cor}
An analytic version of this forward expansion, under
polynomial-growth conditions in a complex cone and with the same
coefficient array expressed through 2-associated Stirling numbers,
was obtained by
\citet[Theorem~4]{jacquet_szpankowski_1999}.
Corollary~\ref{cor:ramanujan-derivative-form}, by contrast, is
obtained under real-variable hypotheses and, among derivative
bounds, requires only a bound on the highest derivative involved
in the chosen truncation.

The appearance of the Mahler polynomials in
Corollary~\ref{cor:ramanujan-derivative-form} has a direct
probabilistic interpretation. If \(\eta(R)\) is a Poisson random
variable with mean \(R\), then
\[
\mathbb E e^{z(\eta(R)-R)}
=
e^{R(e^z-1-z)}
=
\sum_{j=0}^{\infty}
\frac{\rho_j(-R)}{j!}z^j.
\]
Consequently,
\[
\mathbb E\bigl(\eta(R)-R\bigr)^j
=
\rho_j(-R).
\]
Thus, formally Taylor-expanding \(\varphi(\eta(R))\) about its mean
\(R\) suggests
\[
\begin{aligned}
\mathbb E\varphi(\eta(R))
&\approx
\mathbb E\!\left[
\sum_{j=0}^{N}
\frac{\varphi^{(j)}(R)}{j!}
\bigl(\eta(R)-R\bigr)^j
\right]
\\
&=
\sum_{j=0}^{N}
\frac{\varphi^{(j)}(R)}{j!}\rho_j(-R).
\end{aligned}
\]
This is the probabilistic form, made explicit by
\citet{yu_2009}, of the Taylor-polynomial idea underlying Berndt's
proof. Our argument proceeds differently: we first establish the
finite-difference expansion in
Theorem~\ref{thm:ramanujan-real-finite-difference} and then convert it
into the derivative form above, with a rigorous remainder estimate.

To compare Corollary~\ref{cor:ramanujan-derivative-form} with
Ramanujan's expansion~\eqref{eq:ramanujan-expansion}, set \(N=2M-1\).
The terms corresponding to \(j=0\) and \(j=1\) agree trivially, since
\(\rho_0(-R)=1\) and \(\rho_1(-R)=0\). For \(j\geqslant2\), write
\[
\rho_j(-R)
=
\sum_{q=1}^{\lfloor j/2\rfloor}c_{j,q}R^q.
\]
Taking the \(q\)-th derivative at \(x=0\) of both sides of the
generating function~\eqref{eq:mahler-generating-function} gives
\[
\sum_{j=0}^\infty \frac{c_{j,q}}{j!}z^j
=
\frac{(e^z-1-z)^q}{q!}.
\]
Differentiating this generating-function identity for \(c_{j,q}\) leads
to the recurrence relation
\[
c_{j+1,q}
=
q\,c_{j,q}
+
j\,c_{j-1,q-1}.
\]
Consequently, after setting \(k=j-q+1\), the coefficients
\(c_{j,q}\) satisfy the recurrence and initial conditions
in~\eqref{eq:ramanujan-coefficient-recurrence}. Hence
\[
c_{j,q}=b_{j-q+1,j},
\]
or equivalently,
\[
\rho_j(-R)
=
\sum_{k=\lceil(j+2)/2\rceil}^{j}
b_{k,j}R^{j-k+1},
\qquad j\geqslant2.
\]

For a fixed derivative order \(j\), Ramanujan's expansion retains only
the summands with \(k\leqslant M\). Since \(q=j-k+1\), these are
precisely the monomials for which
\[
q\geqslant j-M+1.
\]
On the other hand, the derivative bound~\eqref{eq:ramanujan-derivative-bound} gives
\[
\left|
c_{j,q}R^q\frac{\varphi^{(j)}(R)}{j!}
\right|
\leqslant
|c_{j,q}|C^jG(R)R^{q-j}.
\]
Consequently, every monomial with \(q\leqslant j-M\) contributes
\(O(G(R)R^{-M})\) and may be absorbed into the remainder. Since every
nonzero monomial of \(\rho_j(-R)\), for \(j\geqslant1\), has degree at
least one, no monomial is absorbed when \(j\leqslant M\). Thus the full
coefficient of \(\varphi^{(j)}(R)/j!\) agrees in the two expansions for
\(j\leqslant M\). When \(j>M\), Ramanujan's expansion contains exactly
the monomials of \(\rho_j(-R)\) with
\(q\geqslant j-M+1\), while the remaining lower-degree monomials are of
remainder order. Finally, for \(j=2M-1\),
\[
\deg \rho_j(-R)
\leqslant
\left\lfloor\frac{j}{2}\right\rfloor
=
M-1
=
j-M,
\]
so the entire term involving
\(\varphi^{(2M-1)}(R)\) is of remainder order. The additional
exponentially small remainder in
Corollary~\ref{cor:ramanujan-derivative-form} does not affect this coefficient
comparison.

There are also nontrivial examples of functions \(\varphi(x)\) for
which the asymptotic series of
Corollary~\ref{cor:ramanujan-derivative-form} converges and the
identity
\begin{equation}
\label{eq:ramanujan-convergent-identity}
e^{-R}\sum ^{\infty }_{m=0}\frac{\varphi \left( m\right) }{m!}R^{m}=\sum ^{\infty }_{m=0}\frac{\varphi ^{\left( m\right) }\left( R\right) }{m!}\rho_{m}\left( -R\right) 
\end{equation}
holds for all \(R>0\). More precisely, in Proposition \ref{prop:moment-generating-ramanujan-identity} we show that if \(\varphi\) is the moment-generating transform of a finite measure with a sufficiently light tail, then the above identity holds.

The remainder of the paper is organized as follows. Section~\ref{sec:direct-depoissonization} proves the direct depoissonization estimates announced above, treating first-order estimates before the higher-order Poisson--Charlier expansion. Section~\ref{sec:reverse-depoissonization} establishes the reverse depoissonization theorem for arbitrary positive arguments and then proves the announced Ramanujan-type expansions. Section~\ref{sec:estimating-E} develops sufficient conditions and general tools for estimating the remainder operator \(E\), while Section~\ref{sec:applications} presents applications to summatory sequences, poissonized expectations, Stirling numbers of the second kind, and symmetric tries. Appendix~\ref{app:polynomial-properties} provides a self-contained exposition of the key properties of the special polynomials used throughout the paper, together with proofs of several more specialized facts used in our arguments. Appendix~\ref{app:convergent-expansion-examples} contains nontrivial examples of classes of sequences for which the series \eqref{eq:ramanujan-inverted-series} and \eqref{eq:ramanujan-convergent-identity} converge.
\section{Direct depoissonization estimates}
\label{sec:direct-depoissonization}
This section proves the direct depoissonization results announced in the Introduction. The first-order estimates lead to Corollary~\ref{cor:monotone-depoissonization}, while the higher-order Poisson--Charlier estimates yield Theorem~\ref{thm:main-depoissonization}.

\subsection{First-order estimates}
Throughout this section, unless stated otherwise, the function \(f\) is assumed to be entire.
The sequence \(A_0, A_1, \ldots\) denotes the coefficients in the Taylor expansion of \(e^x f(x)\); equivalently,  
\[
f(x) = e^{-x}\sum_{m=0}^{\infty} \frac{A_m}{m!}\,x^m.
\]
\begin{thm}
\label{thm:first-order-identity}
For any positive \(R\)  the identity
	\[
	\begin{aligned} A_n-f(R)
&=\sum ^{n-1 }_{m=0}\frac{\Delta A_{m}}{m!}\int ^{\infty}_{R}e^{-t}t^{m}\,dt-\sum ^{\infty }_{m= n}\frac{\Delta A_{m}}{m!}\int ^{R}_{0}e^{-t}t^{m}\,dt
\end{aligned}
	\]
	holds for all \(n\geqslant 1\).
\end{thm}
\begin{proof}
	We can express
\[
A_{n}=A_0+\sum ^{n-1}_{m=0}\left( A_{m+1}-A_{m}\right)= A_{0}+\sum ^{n-1}_{m=0}\Delta A_{m} = A_{0}+\sum ^{n-1}_{m=0}\frac{\Delta A_{m}}{m!}m!.
\]
Writing \(m!\) as an Euler integral of the second kind gives the following representation of \(A_n\) in terms of the differences \(\Delta A_m\):
\begin{equation}
\label{eq:coefficient-integral-representation}
A_n = A_{0}+\sum ^{n-1}_{m=0}\frac{\Delta A_{m}}{m!} \int ^{\infty}_{0}e^{-t}t^{m}\,dt.	
\end{equation}
On the other hand
\[
f(R)=f(0)+\int_{0}^Rf'(x)\,dx=A_0+\int_{0}^R\left(e^{-x}\sum ^{\infty }_{m=0}\frac{\Delta A_{m}}{m!}x^{m}
\right)\,dx
\]
hence
\begin{equation}
\label{eq:poisson-transform-integral-representation}
f(R)=A_0+\sum ^{\infty }_{m=0}\frac{\Delta A_{m}}{m!}\int_{0}^R e^{-x}x^{m}
\,dx.	
\end{equation}
Subtracting identity (\ref{eq:poisson-transform-integral-representation}) from (\ref{eq:coefficient-integral-representation}) we obtain the identity of the theorem.
\end{proof}
\begin{thm}
\label{thm:first-order-bound}
	The inequality 
\begin{equation}
\label{eq:first-order-bound-general}
\begin{aligned} \left|A_n-f(R)\right|
&\leqslant
 n!\frac{e^R}{R^n}E(f';R)
\end{aligned}	
\end{equation}
holds for all \(R>0\) and \(n\geqslant 1\). For the special case when \(R=n\) we have 
\begin{equation}
\label{eq:first-order-bound-diagonal}
\begin{aligned} \left|A_n-f(n)\right|
&\leqslant
 2\sqrt{n}E(f';n),
\end{aligned}
\end{equation}
if \(n\geqslant 1\).
\end{thm}
\begin{proof}
From Theorem~\ref{thm:first-order-identity},
\[
\begin{aligned}
|A_n-f(R)|
&\leqslant \sum_{m=0}^{n-1}\frac{|\Delta A_m|}{m!}\int_R^\infty e^{-t}t^m\,dt
   +\sum_{m=n}^{\infty}\frac{|\Delta A_m|}{m!}\int_0^R e^{-t}t^m\,dt .
\end{aligned}
\]
For $R>0$, $t>0$, and integers $m,n\geqslant 0$ we have
\begin{equation}
\label{eq:power-bound-upper-tail}
t^{m} \leqslant t^{m}\Bigl(\frac{t}{R}\Bigr)^{n-m}= R^{\,m-n}t^{n}, \quad \text{if } t\geqslant R,\; m<n,
\end{equation}
and
\begin{equation}
\label{eq:power-bound-lower-tail}
t^{m} \leqslant t^{m}\Bigl(\frac{R}{t}\Bigr)^{m-n}= R^{\,m-n}t^{n}, \quad \text{if } 0<t\leqslant R,\; m\geqslant n.
\end{equation}
Hence
\[
\begin{aligned}
|A_n-f(R)|
&\leqslant \frac{1}{R^n}\sum_{m=0}^{n-1}\frac{|\Delta A_m|}{m!}R^m
            \int_R^\infty e^{-t}t^n\,dt
   +\frac{1}{R^n}\sum_{m=n}^{\infty}\frac{|\Delta A_m|}{m!}R^m
            \int_0^R e^{-t}t^n\,dt .
\end{aligned}
\]
Taking the maximum of the two integrals and recalling the definition \eqref{eq:E-definition} of the operator $E$, we obtain
\[
|A_n-f(R)|
\;\leqslant\; \frac{e^R E(f';R)}{R^n}
\max\!\left\{\int_R^\infty e^{-t}t^n\,dt,\;\int_0^R e^{-t}t^n\,dt\right\}.
\]
Finally, bounding the maximum by the sum, which equals $\Gamma(n+1)=n!$, yields inequality~\eqref{eq:first-order-bound-general}.
For the special case \(R=n\), put
\[
U_n:=\int_n^\infty e^{-t}t^n\,dt,
\qquad
L_n:=\int_0^n e^{-t}t^n\,dt.
\]
Define
\[
H(y):=\frac{\log(1+y)-y}{y^2},
\qquad -1<y<\infty,\quad y\ne0.
\]
Since
\[
H(y)
=
-\frac1{y^2}\int_0^y\frac{t}{1+t}\,dt
=
-\int_0^1\frac{s}{1+ys}\,ds,
\]
the function \(H\) is strictly increasing on \((-1,\infty)\).

The substitutions \(t=n+\sqrt n\,u\) and
\(t=n-\sqrt n\,u\), respectively, give
\[
\begin{aligned}
\frac{1}{\sqrt n}\left(\frac en\right)^n U_n
&=
\int_0^\infty
\exp\!\left(
u^2H\left(\frac{u}{\sqrt n}\right)
\right)\,du,\\
\frac{1}{\sqrt n}\left(\frac en\right)^n L_n
&=
\int_0^{\sqrt n}
\exp\!\left(
u^2H\left(-\frac{u}{\sqrt n}\right)
\right)\,du.
\end{aligned}
\]
For \(0<u<\sqrt n\), the monotonicity of \(H\) gives
\[
H\left(-\frac{u}{\sqrt n}\right)
<
H\left(\frac{u}{\sqrt n}\right),
\]
and hence \(L_n<U_n\). Moreover, since \(u/\sqrt n\leqslant u\),
\[
\begin{aligned}
\frac{1}{\sqrt n}\left(\frac en\right)^n U_n
=\int_0^\infty
\exp\!\left(
u^2H\left(\frac{u}{\sqrt n}\right)
\right)\,du\leqslant
\int_0^\infty \exp\!\bigl(u^2H(u)\bigr)\,du
=
e\int_1^\infty te^{-t}\,dt
=2.
\end{aligned}
\]
The last integral is the case \(n=1\) of the preceding substitution.
Thus
\[
\frac{1}{\sqrt n}\left(\frac en\right)^n
\max\{U_n,L_n\}
\leqslant2,
\]
with equality exactly when \(n=1\). Substitution into the preceding
estimate with \(R=n\) proves
\eqref{eq:first-order-bound-diagonal}.
\end{proof}
\begin{proof}[Proof of Corollary \ref{cor:monotone-depoissonization}]
	If the sequence \((A_j)_{j\geqslant 0}\) is monotone, then \(\Delta A_j\) has the same sign for all \(j\geqslant 0\). Therefore,
	\[
	E(f';R)=e^{-R} \sum_{m=0}^\infty \frac{|\Delta  A_m|}{m!}\, R^m=\left|e^{-R} \sum_{m=0}^\infty \frac{\Delta  A_m}{m!}\, R^m\right|=|f'(R)|
	\]
	for all \(R>0\).
	Inserting this expression into the inequalities of Theorem \ref{thm:first-order-bound} we complete the proof of the corollary.
\end{proof}	

\subsection{Higher-order Poisson--Charlier estimates}
We now develop the higher-order estimates that culminate in the proof of Theorem~\ref{thm:main-depoissonization}.
\begin{lem}
\label{lem:tail-integral-polynomial}
 If $N,m\geqslant 0$ are integers and $x\in\mathbb{R}$, then 
\[
e^{x}\int_x ^{\infty }t^{m}\left( x-t\right) ^{N}e^{-t}\,dt= (-1)^N\sum ^{m}_{j=0}\binom{m}{j}\left( N+m-j\right) !x^{j}.
\]	
\end{lem}
\begin{proof}
With the substitution \(t=x+u\), for every real \(x\) we obtain
\[
\begin{aligned}
e^x\int_x^\infty t^m(x-t)^N e^{-t}\,dt
&=(-1)^N\int_0^\infty (x+u)^m u^N e^{-u}\,du \\
&=(-1)^N\sum_{j=0}^m\binom mj x^j
  \int_0^\infty u^{N+m-j}e^{-u}\,du \\
&=(-1)^N\sum_{j=0}^m\binom mj (N+m-j)!x^j.
\end{aligned}
\]
This is the asserted identity.
\end{proof}
\begin{lem}
\label{lem:coefficient-integral-splitting} If \(m,n,N\geqslant0\) are integers then
\begin{equation}
\label{eq:coefficient-integral-cases}
	[x^n]e^x\int ^{x}_{0} t^{m} e^{-t} \left( x-t\right) ^{N}\,dt=
		\begin{cases}
            \displaystyle [x^n]e^x\int ^{\infty}_{0} t^{m} e^{-t} \left( x-t\right) ^{N}\,dt, & \text{if $m<n$}
            \\
            	0, & \text{if $m\geqslant n$ }
		 \end{cases}
\end{equation}	
\end{lem}
\begin{proof}

\emph{Case $m\geqslant n$.}	Note that making a change of variables \(t\to xt\) in the appropriate integral and expanding the resulting \(e^{x(1-t)}\) exponential term under the integration sign into a Taylor series we get
\[
\begin{aligned}
e^x\int ^{x}_{0} t^{m} e^{-t} \left( x-t\right) ^{N}\,dt
&=x^{m+N+1}\int ^{1}_{0} t^{m} e^{x(1-t)}
  \left( 1-t\right) ^{N}\,dt \\
&=x^{m+N+1}\sum_{j=0}^\infty\frac{x^j}{j!}
  \int ^{1}_{0} t^{m}\left( 1-t\right) ^{N+j}\,dt.
\end{aligned}
\]
The series on the right starts at degree $m+N+1$, 
which means that  coefficients of powers of  \(x\) in the Taylor expansion of this function for powers smaller than \(m+N+1\) are zero.  This proves the statement of the lemma when \(m\geqslant n\).

\emph{Case $m<n$.}
 Lemma \ref{lem:tail-integral-polynomial} implies that
\[
[x^n]e^x\int ^{\infty}_{x} t^{m} e^{-t} \left( x-t\right) ^{N}\,dt=[x^n](-1)^N\sum ^{m}_{j=0}\binom{m}{j}\left( N+m-j\right) !x^{j}=0
\]
whenever \(m<n\). Therefore, when $m<n$, the contribution of the integral over $(x,\infty)$ to the coefficient of $x^n$ vanishes, and we obtain
\[
\begin{split}
[x^n]e^x\int ^{x}_{0} t^{m} e^{-t} \left( x-t\right) ^{N}\,dt&=[x^n]e^x\int ^{x}_{0} t^{m} e^{-t} \left( x-t\right) ^{N}\,dt+[x^n]e^x\int ^{\infty}_{x} t^{m} e^{-t} \left( x-t\right) ^{N}\,dt
\\
&=[x^n]e^x\int ^{\infty}_{0} t^{m} e^{-t} \left( x-t\right) ^{N}\,dt	
\end{split}
\]
as claimed.
\end{proof}
For real \(R\ne0\) and nonnegative integers \(m,n\), we adopt the notation
\begin{equation}
\label{eq:charlier-coefficient-definition}
R^mC_m(R,n)=n!\left[ z^{n}\right] e^{z}\left( z-R\right) ^{m}.	
\end{equation}
For fixed \(R\ne0\) and nonnegative integer \(m\), the functions \(C_m(R,x)\) are polynomials of degree \(m\) in the variable \(x\), known as \emph{Charlier polynomials}. When \(R>0\), these polynomials form a system of orthogonal polynomials with respect to the Poisson distribution with parameter \(R\).
\begin{thm}
\label{thm:charlier-remainder-identity} For every \(R>0\) and all integers \(n\geqslant1\) and \(N\geqslant0\), the following identity
	\[
\begin{aligned}
A_n -\sum ^{N}_{m=0}\frac{f^{\left( m\right) }\left( R\right) }{m!}R^{m}C_{m}(R,n)
&= \frac{1}{N!} \sum ^{n-1 }_{m=0}\frac{\Delta^{N+1}A_{m}}{m!}\int_R^\infty t^m e^{-t} t^{N}C_{N}(t,n)\,dt
\\
&\quad-\frac{1}{N!} \sum ^{\infty }_{m=n  }\frac{\Delta^{N+1}A_{m}}{m!}\int_0^Rt^m e^{-t} t^{N}C_{N}(t,n)\,dt
\end{aligned}
\]
holds.
\end{thm}
\begin{proof}
By Taylor's theorem with the integral form of the remainder, we have
\[
\begin{split}
f(x)-\sum ^{N}_{m=0}\frac{f^{\left( m\right) }\left( R\right) }{m!}\left( x-R\right) ^{m}=\frac{1}{N!}\int ^{x}_{R}f^{\left( N+1\right) }\left( t\right) \left( x-t\right) ^{N}\,dt.	
\end{split}
\]
Multiplying both sides of this identity by \(e^x\), expressing the integral on the right as the difference of two integrals (over \([0,x]\) and \([0,R]\))  and  extracting the coefficient at \(x^n\) in the expansion of both sides we get
\[
\begin{aligned}
&A_n -\sum ^{N}_{m=0}\frac{f^{\left( m\right) }\left( R\right) }{m!}n!\left[ x^{n}\right] e^{x}\left( x-R\right) ^{m}
\\
&\quad= \frac{n!}{N!} [x^n]e^x\int ^{x}_{0}f^{\left( N+1\right) }\left( t\right) \left( x-t\right) ^{N}\,dt-\frac{n!}{N!}\left[ x^{n}\right]e^{x}\int_0^R{f^{\left( N+1\right) }\left( t\right) } \left( x-t\right) ^{N}\,dt.
\end{aligned}
\]
Expanding \(f^{(N+1)}(t)\) into a series
\[
f^{(N+1)}(t)=e^{-t}\sum ^{\infty }_{m=0}\frac{\Delta^{N+1}A_{m}}{m!}t^m
\]
and exchanging summation and integration signs we get
\begin{equation*}
\begin{aligned}
&A_n -\sum ^{N}_{m=0}\frac{f^{\left( m\right) }\left( R\right) }{m!}n!\left[ x^{n}\right] e^{x}\left( x-R\right) ^{m}
\\
&\quad= \frac{n!}{N!} \sum ^{\infty }_{m=0}\frac{\Delta^{N+1}A_{m}}{m!}[x^n] e^x\int ^{x}_{0} t^{m} e^{-t} \left( x-t\right) ^{N}\,dt
\\
&\qquad-\frac{n!}{N!} \sum ^{\infty }_{m=0}\frac{\Delta^{N+1}A_{m}}{m!} \left[ x^{n}\right]e^x\int_0^Rt^m e^{-t}\left( x-t\right) ^{N}\,dt.
\end{aligned}
\end{equation*}
By Lemma \ref{lem:coefficient-integral-splitting}, all terms with \(m\geqslant n\) vanish in the first sum and  the upper limits of integration of the integrals inside the remaining part of the sum can be extended from \(x\) to \(+\infty\) according to the same lemma, leading to the expression
\[
\sum ^{\infty}_{m=0}\frac{\Delta^{N+1}A_{m}}{m!}[x^n] e^x\int ^{x}_{0} t^{m} e^{-t} \left( x-t\right) ^{N}\,dt=\sum ^{n-1 }_{m=0}\frac{\Delta^{N+1}A_{m}}{m!}[x^n] e^x\int ^{\infty}_{0} t^{m} e^{-t} \left( x-t\right) ^{N}\,dt.
\]
This, after some straightforward cancellations of integral terms and recalling the definition \eqref{eq:charlier-coefficient-definition} of Charlier polynomials, leads to the identity of the theorem.
\end{proof}

As noted in the Introduction, the diagonal values of the Charlier
polynomials satisfy
\[
\tau_m(n)=n^m C_m(n,n),
\qquad m\geqslant0,\quad n\geqslant1.
\]
Conversely, Charlier polynomials can be expressed in terms of the
polynomials \(\tau_m\), as shown in
Lemma~\ref{lem:charlier-tau-duality}.
Proofs of Lemmas~\ref{lem:charlier-pointwise-bound} and~\ref{lem:charlier-integral-bound}, together with a self-contained exposition of the properties of Charlier polynomials needed in this paper, are given in Appendix~\ref{app:polynomial-properties}. The reverse identities in Lemmas~\ref{lem:charlier-inversion}, \ref{lem:charlier-tau-shift}, and~\ref{lem:mahler-tau-stirling} are stated in Section~\ref{sec:reverse-depoissonization}; their proofs are also given in Appendix~\ref{app:polynomial-properties}.
For further background and additional references on these polynomials, see, for example, \citet{szego_1975} and \citet{chihara_1978}.
\begin{lem}
\label{lem:charlier-pointwise-bound} If \(|y|\geqslant 1\) then for all real \(t\) and integers \(N\geqslant 0\) one has
\[
\bigl|t^NC_N(t,y)\bigr|\leqslant  N!|y|^{N/2}e^{\frac{|y-t|}{\sqrt{|y|}}}.
\]
At \(t=0\), the product \(t^NC_N(t,y)\) is understood by polynomial continuation from \(t\ne0\).
\end{lem}

\par\smallskip
\begin{lem}
\label{lem:charlier-integral-bound}
If \(n\geqslant 2\) and \( N\geqslant 0\) then
\[
\int_0^\infty t^{n} e^{-t}\bigl|t^NC_N(t,n)\bigr|\,dt \leqslant 12N!n!n^{N/2}.
\]
\end{lem}
\begin{thm}
	\label{thm:charlier-depoissonization-bound}
 For every \(R>0\) and all integers \(n\geqslant 2\) and \(N\geqslant 0\), the following inequality holds	\[
	\begin{aligned}
\left|A_n -\sum ^{N}_{m=0}\frac{f^{\left( m\right) }\left( R\right) }{m!}R^mC_m(R,n)\right|
\leqslant 12 n!\frac{e^RE(f^{(N+1)};R)}{R^{n}}n^{N/2}.
\end{aligned}
	\]
\end{thm}
\begin{proof} Starting from the identity in Theorem \ref{thm:charlier-remainder-identity}, and evaluating the term \(t^m\) under the integration sign by means of inequalities \eqref{eq:power-bound-upper-tail} and \eqref{eq:power-bound-lower-tail} we obtain
	\[
	\begin{aligned}
&\left|A_n -\sum ^{N}_{m=0}\frac{f^{\left( m\right) }\left( R\right) }{m!}R^mC_m(R,n)\right|
\\
&\leqslant \frac{1}{N!} \sum ^{n-1 }_{m=0}\frac{|\Delta^{N+1}A_{m}|}{m!}\int_R^\infty t^m e^{-t} \left(\frac{R}{t}\right)^{m-n}\bigl|t^NC_N(t,n)\bigr|\,dt
\\ &\quad+\frac{1}{N!} \sum ^{\infty }_{m=n}\frac{|\Delta^{N+1}A_{m}|}{m!}\int_0^Rt^m e^{-t}  \left(\frac{R}{t}\right)^{m-n}\bigl|t^NC_N(t,n)\bigr|\,dt
\\
&\leqslant \frac{e^RE(f^{(N+1)};R)}{N!R^{n}} \int_0^\infty t^{n} e^{-t}\bigl|t^NC_N(t,n)\bigr|\,dt.
\end{aligned}
	\]
	Applying Lemma \ref{lem:charlier-integral-bound} to evaluate the integral in the last expression completes the proof.
\end{proof}
\begin{proof}[Proof of Theorem \ref{thm:main-depoissonization}]
Substituting \(R=n\) in the inequality of Theorem~\ref{thm:charlier-depoissonization-bound} and using the Stirling–type bound
\[
n! \leqslant e\sqrt{n}\left(\frac{n}{e}\right)^n, \qquad n\geqslant 1,
\]
we obtain the stated estimate. Finally, the numerical factor is simplified using 
\(12e < 33\), which completes the proof.
\end{proof}
\section{Reverse depoissonization and\texorpdfstring{\\}{ }Ramanujan-type expansions}
\label{sec:reverse-depoissonization}
This section proves Theorems~\ref{thm:ramanujan-integer-finite-difference} and~\ref{thm:ramanujan-real-finite-difference} and Corollary~\ref{cor:ramanujan-derivative-form}, which were announced in the Introduction. We first establish the arbitrary-real-argument result, Theorem~\ref{thm:reverse-depoissonization-arbitrary-R}; the integer case then follows by setting \(R=n\).

The inversion step uses the following identity.
\begin{lem}
\label{lem:charlier-inversion} For every real \(R\ne0\), every real \(n\), and every integer \(s \geqslant 0\),
	\[
\sum ^{s}_{j=0}\frac{R^{s-j}C_{s-j}(R,n)(-R)^jC_{j}(-R,-n) }{(s-j)!j!}=\begin{cases}
	0, & \text{if $s\geqslant 1$ }
	\\
	1&\text{if $s=0$ }.
\end{cases}
\]
\end{lem}

\begin{thm}
\label{thm:reverse-depoissonization-arbitrary-R}
Let $(A_m)_{m\geqslant 0}$ be a sequence of real numbers such that the
exponential generating series
\[
\sum_{m=0}^{\infty} A_m\frac{x^m}{m!}
\]
converges for every $x\in\mathbb{R}$, and define its Poisson transform
by
\[
f(x):=e^{-x}\sum_{m=0}^{\infty} A_m\frac{x^m}{m!}.
\]
Then, for every integer $N\geqslant 0$, every $R>0$, and every integer
$n\geqslant 2$,
\[
\begin{aligned}
&\left|
f(R)-\sum_{m=0}^{N}\frac{\Delta^m A_n}{m!}
(-R)^m C_m(-R,-n)
\right| \\
&\quad\leqslant
12(N+1)n!\frac{e^R}{R^n}
e^{\frac{|R-n|}{\sqrt{n}}}  n^{N/2}E\!\left(f^{(N+1)};R\right),
\end{aligned}
\]
where
\[
E\!\left(f^{(N+1)};R\right)
:=
e^{-R}\sum_{j=0}^{\infty}
\frac{\left|\Delta^{N+1}A_j\right|}{j!}R^j.
\]
\end{thm}
\begin{proof} Applying Theorem \ref{thm:charlier-depoissonization-bound} to functions \(f(x),f'(x),\ldots,f^{(N)}(x)\) we can rewrite the resulting \(N+1\) inequalities in the following way
	\[
	\begin{aligned}
\Delta^j A_n =\sum ^{N-j}_{m=0}\frac{f^{\left( m+j\right) }\left( R\right) }{m!}R^mC_m(R,n)
+c_j\cdot12n!\frac{e^RE(f^{(N+1)};R)}{R^{n}}n^{(N-j)/2},
\end{aligned}
	\]
	with  \(|c_j|\leqslant 1\).
Multiplying both sides of this identity by \((-R)^jC_j(-R,-n)/j!\), summing over \(j=0,1,2,\ldots,N\) and exchanging the order of summation, we get
\[
\begin{aligned}
\sum_{j=0}^N\frac{\Delta^j A_n}{j!}(-R)^jC_j(-R,-n)
&=\sum_{s=0}^Nf^{(s)}(R) \sum ^{s}_{j=0}
\frac{R^{s-j}C_{s-j}(R,n)(-R)^jC_{j}(-R,-n)}{(s-j)!j!} \\
&\quad+12n!\frac{e^R}{R^{n}}E(f^{(N+1)};R) \\
&\qquad{}\times\sum_{j=0}^Nc_jn^{(N-j)/2}
\frac{(-R)^jC_j(-R,-n)}{j!}.
\end{aligned}
\]
By Lemma \ref{lem:charlier-inversion}, the inner sum equals 1 when \(s=0\) and vanishes otherwise. Hence
\[
\begin{aligned}
\sum_{j=0}^N\frac{\Delta^j A_n}{j!}(-R)^jC_j(-R,-n)
&=f(R) +12n!\frac{e^R}{R^{n}}E(f^{(N+1)};R) \\
&\qquad{}\times\sum_{j=0}^Nc_jn^{(N-j)/2}
\frac{(-R)^jC_j(-R,-n)}{j!}.
\end{aligned}
\]
Applying the upper bound \(\bigl|(-R)^jC_j(-R,-n)\bigr|/j!\leqslant n^{j/2}e^{\frac{|R-n|}{\sqrt{n}}}\) from Lemma \ref{lem:charlier-pointwise-bound} to the error term, we immediately obtain the estimate of the theorem.
\end{proof}

\begin{proof}[Proof of Theorem \ref{thm:ramanujan-integer-finite-difference}]
Apply Theorem~\ref{thm:reverse-depoissonization-arbitrary-R} with \(R=n\), and then use the standard bound
\[
n! \leqslant e\sqrt{n}\left(\frac{n}{e}\right)^n, \qquad n\geqslant 1.
\]
This yields the desired estimate, after simplifying the numerical factor using
\(12e < 33\). The proof is complete.
\end{proof}

To pass from finite differences at \(n=\lfloor R\rfloor\) to finite differences at \(R\) in the proof of Theorem~\ref{thm:ramanujan-real-finite-difference}, we use the following identity.
\begin{lem}
\label{lem:charlier-tau-shift} For every real \(R\ne0\), every real \(n\), and every integer \(j \geqslant 0\), the following identity holds:
	\[
	\sum ^{j}_{s=0}\frac{\left( -R\right) ^{s}C_{s}\left( -R,-n\right) }{s!}\Delta_x ^{s}x^{j} \Big|_{x=0} =\sum ^{j}_{s=0}\frac{\tau_{s}\left( -R\right) }{s!}\Delta_x ^{s}x^{j} \Big|_{x=R-n}.
	\]
\end{lem}

\begin{proof}[Proof of Theorem \ref{thm:ramanujan-real-finite-difference}]Let
\[
A_m:=\varphi(m),\qquad
f(R):=e^{-R}\sum_{m=0}^{\infty}\frac{\varphi(m)}{m!}R^m .
\]
Under the hypotheses of the theorem, $N$ is fixed, and the constants
in the $O(\cdot)$-terms may depend on $N$, $B$, $C$, $T$, and on the
values of $\varphi$ on a fixed bounded interval. We first verify the
growth condition needed below. Fixing \(x=T\) in~\eqref{eq:variation-condition}
shows that, for \(y\geqslant T\),
\[
G(y)\leqslant B G(T)e^{C(y-T)/\sqrt T}\leqslant C_1e^{ay}
\]
for suitable constants \(C_1,a>0\). Hence~\eqref{eq:highest-derivative-bound},
together with continuity on bounded intervals, gives
\(\lvert\varphi^{(N+1)}(y)\rvert\leqslant C_2e^{ay}\) for
\(y\geqslant0\). Taylor's formula with integral remainder then yields
\(\lvert\varphi(y)\rvert\leqslant C_3e^{a'y}\) for some \(a'>0\).
Consequently,
\[
\sum_{m=0}^\infty\frac{|\varphi(m)|}{m!}|z|^m
\leqslant C_3\exp\!\left(e^{a'}|z|\right)<\infty,
\qquad z\in\mathbb C,
\]
so the exponential generating series of \(A_m=\varphi(m)\) has
infinite radius of convergence. Applying
Theorem~\ref{thm:reverse-depoissonization-arbitrary-R} with \(n=\lfloor R\rfloor\)  and \(A_m=\varphi(m)\) provides the estimate
\begin{equation}
\label{eq:reverse-depoissonization-intermediate}
\begin{aligned}
f(R) =\sum ^{N}_{m=0}\frac{\Delta ^{m}\varphi(n) }{m!}(-R)^mC_m(-R,-n)+
O\bigl(R^{(N+1)/2}E(f^{(N+1)};R)\bigr)
\end{aligned}
\end{equation}
where
\[
f(R)=e^{-R}\sum_{m=0}^\infty\frac{\varphi(m)}{m!}R^m.
\]
Let us first evaluate the error term in this estimate. For this aim we will need upper bounds for finite differences of \(A_m=\varphi(m)\). By the conditions \eqref{eq:highest-derivative-bound} and \eqref{eq:variation-condition} of our theorem we have
\[
\begin{split}
|\Delta ^{N+1}\varphi(m)|&\leqslant \max_{t\in[m,m+N+1]}|\varphi^{(N+1)}(t)|
\\
&\leqslant \max_{t\in[m,m+N+1]}\frac{G(t)}{t^{N+1}}
\\
&\leqslant B\max_{t\in[m,m+N+1]}\frac{e^{C\frac{|t-R|}{\sqrt{R}}} }{t^{N+1}}G(R)	\\
&\leqslant B\frac{e^{C\frac{N+1}{\sqrt{R}}}}{m^{N+1}}e^{C\frac{|m-R|}{\sqrt{R}}}G(R)
\end{split}
\]
for all \(m\geqslant T\).
Applying the above estimate to evaluate the coefficients of the series defining \(E(f^{(N+1)};R)\) we get
\[
\begin{split}
E(f^{(N+1)};R)&=e^{-R} \sum ^{\infty }_{m=0}\frac{|\Delta ^{N+1}\varphi(m)|}{m!}R^{m}
\\
&\leqslant C_1e^{-R} \left(\sum_{m<T}\frac{|\Delta ^{N+1}\varphi(m)|}{m!}R^{m}+\frac{G(R)}{R^{N+1}}\sum_{m\geqslant T}{
e^{C\frac{|m-R|}{\sqrt{R}}}}\frac{
R^{m+N+1}}{(m+N+1)!}\right)
\\
&=O(R^Te^{-R})+O\left(\frac{G(R)}{R^{N+1}}\mathbb{E}e^{C\frac{|\eta(R)-R|}{\sqrt{R}}}\right)\end{split}
\]
where \(\eta(R)\) is a Poisson random variable with parameter \(R\) and \(C_1\) is some constant that is bounded when \(N\) is fixed and \(R\to\infty\). Since for the Poisson random variable \(\eta(R)\) the expectation inside the \(O(\ldots)\) is finite,
\[
\mathbb{E}e^{C\frac{|\eta(R)-R|}{\sqrt{R}}}\leqslant\mathbb{E}e^{C\frac{\eta(R)-R}{\sqrt{R}}}+\mathbb{E}e^{-C\frac{\eta(R)-R}{\sqrt{R}}}=e^{R\left( e^{\frac{C}{\sqrt{R}}}-1-\frac{C}{\sqrt{R}}\right) }+e^{R\left( e^{\frac{-C}{\sqrt{R}}}-1-\frac{-C}{\sqrt{R}}\right) }=O(1),
\]
for fixed $C>0$ as $R\to\infty$, we conclude that
\[
\begin{split}
E(f^{(N+1)};R)&=O(R^Te^{-R})+O\left(\frac{G(R)}{R^{N+1}}\right).
\end{split}
\]
Inserting this estimate into the error term of \eqref{eq:reverse-depoissonization-intermediate} we get
\begin{equation}
\label{eq:reverse-depoissonization-at-floor-R}
	f\left( R\right)
	=\sum ^{N}_{m=0}\frac{\Delta^{m}\varphi(n)}{m!}(-R)^mC_m(-R,-n)+O\left(R^{-\frac{N+1}{2}}G(R)\right)+O(R^{T+(N+1)/2}e^{-R}).
\end{equation}
Let us now turn to evaluating \(\Delta^{m}\varphi(n)\)  in terms of \(\Delta^{m}\varphi(R)\) inside the above estimate.
By Taylor's formula with integral form of remainder 
\[
\varphi \left( x\right) =\sum ^{N}_{j=0}\dfrac{\varphi^{\left( j\right) }\left( n\right) }{j!}\left( x-n\right) ^{j}+ \frac{1}{N!}\int ^{x}_{n}\varphi^{\left( N+1\right) }\left( t\right) \left( x-t\right) ^{N}\,dt
\]
evaluating the remainder using the upper bound \eqref{eq:highest-derivative-bound} for the \(N+1\)-th derivative of \(\varphi\) we get an estimate
\[
\varphi \left( x\right) =\sum ^{N}_{j=0}\frac{\varphi^{\left( j\right) }\left( n\right) }{j!}\left( x-n\right) ^{j}+ O\left(\frac{G(R)}{R^{N+1}}\right)\]
 that will be true if \(x\in[n,n+N+1]\). Applying the operator \(\Delta_x^s\) with \(s\leqslant N \) to both sides of this estimate and choosing \(x=n\) and \(x=R\)  we get
 \begin{equation}
 \label{eq:finite-difference-at-R}
 \Delta^s\varphi \left( R\right)  =\sum ^{N}_{j=0}\frac{\varphi^{\left( j\right) }\left( n\right) }{j!} \Delta_x^sx^{j} \Big|_{x=R-n}+ O\left(\frac{G(R)}{R^{N+1}}\right)
 \end{equation}
 and 
 \begin{equation}
 \label{eq:finite-difference-at-floor-R}
\Delta^s\varphi \left( n\right)  =\sum ^{N}_{j=0}\frac{\varphi^{\left( j\right) }\left( n\right) }{j!} \Delta_x^sx^{j} \Big|_{x=0}+ O\left(\frac{G(R)}{R^{N+1}}\right). 	
 \end{equation}
Replacing the finite differences  \(\Delta^s\varphi \left( n\right)\) inside the estimate (\ref{eq:reverse-depoissonization-at-floor-R}) by their expression as sums of derivatives \(\varphi^{(j)}(n)\) as provided in the above identity \eqref{eq:finite-difference-at-floor-R}  we can continue as
 \[
 \begin{aligned}
f(R)
 &= \sum ^{N}_{s=0}\frac{\left( -R\right) ^{s}
 C_{s}\left( -R,-n\right) }{s!} \\
 &\quad{}\times\left(\sum ^{N}_{j=0}
 \frac{\varphi^{\left( j\right) }\left( n\right) }{j!}
 \Delta_x^sx^{j} \Big|_{x=0}
 + O\left(\frac{G(R)}{R^{N+1}}\right)\right) \\
 &\quad+O\left(R^{-\frac{N+1}{2}}G(R)\right)
 +O(R^{T+(N+1)/2}e^{-R}) \\
 &= \sum ^{N}_{j=0}\frac{\varphi^{\left( j\right) }
 \left( n\right) }{j!}
 \sum ^{j}_{s=0}\frac{\left( -R\right) ^{s}
 C_{s}\left( -R,-n\right) }{s!}
 \Delta_x^sx^{j} \Big|_{x=0} \\
 &\quad+O\left(R^{-\frac{N+1}{2}}G(R)\right)
 +O(R^{T+(N+1)/2}e^{-R})
\end{aligned}
 \]
 Here, in the last step, we applied Lemma \ref{lem:charlier-pointwise-bound} to evaluate \(\left( -R\right) ^{s}C_{s}\left( -R,-n\right)=O(s!R^{s/2})\) when  \(|R-n|\leqslant 1\) and \(s\) is bounded from above. Noticing that the inner sum of the last estimate is exactly the left side of the identity of  Lemma \ref{lem:charlier-tau-shift} and replacing it with the right side of the identity we get
 \[
 \begin{aligned}
 f(R)&= \sum ^{N}_{j=0}\frac{\varphi^{\left( j\right) }\left( n\right) }{j!} \sum ^{j}_{s=0}\frac{\tau_{s}\left( -R\right) }{s!}\Delta_x ^{s}x^{j} \Big|_{x=R-n}+O\left(R^{-\frac{N+1}{2}}G(R)\right)+O(R^{T+(N+1)/2}e^{-R})
\\
&=\sum ^{N}_{s=0}\frac{\tau_{s}\left( -R\right) }{s!} \sum ^{N}_{j=0}\frac{\varphi^{\left( j\right) }\left( n\right) }{j!} \Delta_x ^{s}x^{j} \Big|_{x=R-n}+O\left(R^{-\frac{N+1}{2}}G(R)\right)+O(R^{T+(N+1)/2}e^{-R})
\\
&=\sum ^{N}_{s=0}\frac{\tau_{s}\left( -R\right) }{s!} \left( \Delta^s\varphi \left( R\right)  + O\left(\frac{G(R)}{R^{N+1}}\right)\right)+O\left(R^{-\frac{N+1}{2}}G(R)\right)+O(R^{T+(N+1)/2}e^{-R}).
\end{aligned}
 \]
 In the last step we replaced the inner sum involving derivatives of \(\varphi\) by
 \(
 \Delta^s\varphi \left( R\right)  \) according to (\ref{eq:finite-difference-at-R}).
 Finally, Proposition~\ref{prop:tau-bound} gives
 \(|\tau_s(-R)|\leqslant s!R^{s/2}\) for \(R\geqslant1\).
 Since \(s\leqslant N\) and \(N\) is fixed, applying this bound to the
 accumulated \(O(G(R)/R^{N+1})\) terms shows that they are absorbed by
 the stated \(O(R^{-(N+1)/2}G(R))\) remainder. This completes the proof.
 \end{proof}

 To convert the finite-difference expansion into the derivative form of Corollary~\ref{cor:ramanujan-derivative-form}, we use the following relation between Mahler polynomials, the polynomials \(\tau_s\), and Stirling numbers. Its proof is given in
Appendix~\ref{app:polynomial-properties}.
 \begin{lem}
\label{lem:mahler-tau-stirling}
For all \(x \in \mathbb{R}\) and all integers \(j \geqslant 0\), Mahler's polynomials \(\rho_j(x)\) defined by the identity \eqref{eq:mahler-generating-function} can be expressed in terms of polynomials \(\tau_s(x)\) as
	 \[
 \rho_j(x)=\sum ^{j }_{s=0}\tau_s(x)\Stirling{j}{s},
 \]
where \(\Stirling{j}{s}\) denotes the Stirling number of the second kind.
\end{lem}

 \begin{proof}[Proof of Corollary \ref{cor:ramanujan-derivative-form}]
	By Taylor's formula with the integral form of the remainder, uniformly for \(x\in[R,R+N]\),
 \[
\varphi \left( x\right) =\sum ^{N}_{j=0}\frac{\varphi^{\left( j\right) }\left( R\right) }{j!}\left( x-R\right) ^{j}+ O\left(\frac{G(R)}{R^{N+1}}\right).
\]
Applying the finite-difference operator \(\Delta_x^s\) with \(s\leqslant N\) to both sides and evaluating at \(x=R\) gives
\[
\Delta^s\varphi \left( R\right) =\sum ^{N}_{j=s}\frac{\varphi^{\left( j\right) }\left( R\right) }{j!} \Delta_x^sx^{j} \Big|_{x=0}+O\left(\frac{G(R)}{R^{N+1}}\right).
\]
Using the standard identity for finite differences in terms of Stirling numbers of the second kind (see [1b], p.~204 of \citet{comtet_1974}),
\[
\Delta^s_xx^{j} \Big|_{x=0}=s!\Stirling{j}{s}
\]
we obtain
\[
\Delta^s\varphi \left( R\right) =s!\sum ^{N}_{j=s}\frac{\varphi^{\left( j\right) }\left( R\right) }{j!} \Stirling{j}{s}+O\left(\frac{G(R)}{R^{N+1}}\right).
\]
Let us insert this into the main term of Theorem~\ref{thm:ramanujan-real-finite-difference} and exchange the order of summation:
 \begin{equation}
 \label{eq:ramanujan-derivative-intermediate}
 \begin{aligned}
 e^{-R}\sum_{m=0}^\infty\frac{\varphi(m)}{m!}R^m
 &=\sum ^{N}_{j=0}\frac{\varphi ^{\left( j\right) }
 \left( R\right) }{j!}
 \sum ^{j}_{s=0}{\tau_{s}\left( -R\right) }\Stirling{j}{s} \\
 &\quad+O\left(R^{-\frac{N+1}{2}}G(R)\right)
 +O(R^{T+(N+1)/2}e^{-R}).
 \end{aligned}
 \end{equation}
According to Lemma \ref{lem:mahler-tau-stirling} we can replace here the sum
 \[
 \sum ^{j}_{s=0}{\tau_{s}\left( -R\right) }\Stirling{j}{s}
 \]
with \(\rho_j(-R)\), which completes the proof of the corollary. The implied constants have the dependencies specified in Theorem~\ref{thm:ramanujan-real-finite-difference}.
 \end{proof}
\section{Estimating the remainder operator \texorpdfstring{\(E\)}{E}}
\label{sec:estimating-E}
The direct and reverse depoissonization estimates reduce applications to obtaining effective bounds for \(E(f^{(N+1)};R)\). We first record two elementary sufficient conditions and then develop tools for functions satisfying functional or differential equations.

\subsection{Elementary sufficient conditions}

At first glance, since under our setting the behavior of the coefficients \(A_m\) is assumed to be unknown while only information about the analytic behavior of \(f(x)\) is available, it may seem difficult to obtain upper bounds for the finite differences \(|\Delta A_m|\), and consequently for \(E(f^{(N+1)};R)\). 
Although this is generally true, there are two important cases in which such bounds are in fact easy to obtain.

\begin{enumerate}
\item
For a fixed integer \(N\geqslant0\), suppose that either
\(
\Delta^{N+1}A_{m}\geqslant 0
\)
for all \(m\geqslant 0\), or 
\(
\Delta^{N+1}A_{m}\leqslant 0
\)
for all \(m\geqslant 0\). 
In this case, for every \(R\geqslant0\), we have \(E(f^{(N+1)};R)=\bigl|f^{(N+1)}(R)\bigr|\).
Since information about the asymptotic behavior of \(f(x)\) and its derivatives is typically available, obtaining an upper bound for \(E(f^{(N+1)};R)\) in this situation is straightforward.

\item 
Suppose that for some nonnegative integers \(s\) and \(k\) we can represent
\[
\Delta^s A_m = h(m),
\]
where \(h(x)\) is a function that is continuous and possesses \(k+1\) continuous derivatives on \([0,\infty)\). Then
\[
|\Delta^{s+k+1} A_m| = |\Delta^{k+1} h(m)| 
   \leqslant \max_{x\in[m,m+k+1]} |h^{(k+1)}(x)|.
\]
Consequently, for every \(R\geqslant0\),
\[
E(f^{(s+k+1)};R)
   \leqslant e^{-R}\sum_{m=0}^\infty
   \frac{R^m}{m!}\,
   \max_{x\in[m,m+k+1]} |h^{(k+1)}(x)|.
\]
Whenever the series on the right converges, it gives a finite upper
bound for \(E(f^{(s+k+1)};R)\). When the derivatives of \(h\) admit
suitable estimates, this majorant can be bounded effectively.
\end{enumerate}

\subsection{Algebraic properties of \texorpdfstring{\(E\)}{E}}
Another strategy for estimating \(E(f^{(N+1)};R)\) applies when the function \(f\) satisfies a functional or differential equation.  
In that case, one can apply the operator \(E\) to the equation itself and use the algebraic properties described in Theorem~\ref{thm:E-properties}.
\begin{thm}
\label{thm:E-properties}
Let \(f\) and \(g\) be entire functions, and let \(p,q \geqslant 0\) with \(p+q = 1\). Then, the first three inequalities below hold for all \(r \geqslant 0\), while the fourth holds for \(r>0\):
\begin{enumerate}
\item For any \(a,b \in \mathbb{R}\),
    \begin{equation} \label{eq:E-linear-combination}
        E_x\bigl(a f(x) + b g(x); \, r\bigr) \leqslant |a| \, E(f; r) + |b| \, E(g; r).
    \end{equation}
    \item 
    \begin{equation} \label{eq:E-product}
        E_x\bigl(f(px) \, g(qx); \, r\bigr) \leqslant E(f; rp) \, E(g; rq).
    \end{equation}
    \item
    \begin{equation} \label{eq:E-scaling}
        E_x\bigl(f(px); \, r\bigr) \leqslant E(f; rp).
    \end{equation}
 \item
 \begin{equation} \label{eq:E-derivative-bound}
\left|\frac{d}{dr}E(f;r)\right|\leqslant E(f';r).
 \end{equation}
\end{enumerate}
\end{thm}

\begin{proof}The inequality (\ref{eq:E-linear-combination}) is trivial.
If
\[
f\left( x\right) =e^{-x}\sum ^{\infty }_{m=0}\frac{A_{m}}{m!}x^{m}\quad\hbox{and}\quad g\left( x\right) =e^{-x}\sum ^{\infty }_{m=0}\frac{B_{m}}{m!}x^{m}
\]
then
\[\begin{split}
f(xp)g(xq)&=e^{-xp}e^{-xq}\sum ^{\infty }_{j=0}\frac{A_{j}}{j!}(xp)^{j}\sum ^{\infty }_{m=0}\frac{B_{m}}{m!}(xq)^{m}
\\
&=e^{-x}\sum_{n=0}^\infty\left(\sum_{j=0}^n\frac{n!}{j!(n-j)!}A_jp^jB_{n-j}q^{n-j}\right)\frac{x^n}{n!}
\end{split}
\]
hence
\[
\begin{split}
E_x\bigl(f(xp)g(xq);r\bigr)&=e^{-r}\sum_{n=0}^\infty\left|\sum_{j=0}^n\frac{n!}{j!(n-j)!}A_jp^jB_{n-j}q^{n-j}\right|\frac{r^n}{n!}
\\
&\leqslant e^{-r}\sum_{n=0}^\infty\frac{r^n}{n!}\sum_{j=0}^n\frac{n!}{j!(n-j)!}|A_j|p^j|B_{n-j}|q^{n-j}
\\
&=e^{-rp}\sum ^{\infty }_{m=0}\frac{|A_{m}|}{m!}(rp)^{m}e^{-rq}\sum ^{\infty }_{m=0}\frac{|B_{m}|}{m!}(rq)^{m}
\\
&=E(f;rp)E(g;rq).
\end{split}
\]
This establishes the inequality (\ref{eq:E-product}). The inequality (\ref{eq:E-scaling}) follows by choosing \(g(x)\equiv 1\) and noticing that \(E(1;r)\equiv 1\).  The inequality \eqref{eq:E-derivative-bound} follows from the trivial inequality
\[
\bigl|\Delta |A_m|\bigr|\leqslant |\Delta A_m|.
\]
\end{proof}

\subsection{Bounds derived from functional equations}
As an illustration of how Theorem~\ref{thm:E-properties} can be applied, we now establish Theorem~\ref{thm:functional-equation-bound}.
\begin{thm}
\label{thm:functional-equation-bound}
Fix an integer \(N\geqslant0\). Let \(f\) be an entire function and define
\[
g(z):=f(z)-2f\left(\frac z2\right).
\]
Assume that
	\begin{equation}
	\label{eq:functional-equation-E-condition}
	E(g^{(N+1)};r)=O\left(\frac{1}{r^N}\right)	
	\end{equation}
	as \(r\to +\infty\). Then
	\[
	A_n -\sum ^{N}_{m=0}\frac{f^{\left( m\right) }\left( n\right) }{m!}\tau_m(n)=O\left(\frac{\log n}{n^{\frac{N-1}{2}}}\right),
	\]
	as \(n\to +\infty\).
The implied asymptotic constant may depend on the fixed integer \(N\).
\end{thm}
\begin{proof}
By the definition of \(g\),
\[
f(z)=2f\left(\frac z2\right)+g(z).
\]
Differentiating this identity \(N+1\) times gives
\[
f^{\left( N+1\right) }\left( x\right) =\dfrac{1}{2^{N}}f^{\left( N+1\right) }\left( \dfrac{x}{2}\right) +g^{\left( N+1\right) }\left( x\right) 
\]
Applying the operator \(E_x\) to both sides of this equation and using inequalities (\ref{eq:E-linear-combination})  and (\ref{eq:E-scaling}) to evaluate the right side we get
\[
E\left(f^{\left( N+1\right) };r\right)\leqslant\frac{1}{2^{N}}E\left(f^{\left( N+1\right) };\frac r2\right)+E\left(g^{\left( N+1\right) };r\right).
\]
Iterating this inequality \(M\) times we get
\[
E\left(f^{\left( N+1\right) };r\right)\leqslant\dfrac{1}{2^{N(M+1)}}E\left(f^{\left( N+1\right) }; \frac{r}{2^{M+1}}\right)+\sum_{j=0}^{M}\frac{1}{2^{jN}}E\left(g^{\left( N+1\right) };\frac r{2^j}\right).
\]
First suppose that \(N=0\). For \(r\geqslant2\), choose \(M=\lfloor\log_2 r\rfloor\) in the last inequality. Then the argument of the remaining term belongs to \([1/2,1)\), so that term is bounded. Moreover, the hypothesis gives \(E(g';u)=O(1)\) for \(u\geqslant1\), after enlarging the bound on a compact interval. Since the sum contains \(M+1=O(\log r)\) terms, it follows that
\[
E(f';r)=O(\log r).
\]
This is the required estimate when \(N=0\).

Now suppose that \(N\geqslant1\). Letting \(M\to \infty\) and taking into account that the first summand on the right goes to zero, we get
\[
E\left(f^{\left( N+1\right) };r\right)\leqslant\sum_{j=0}^\infty\frac{1}{2^{jN}}E\left(g^{\left( N+1\right) };\frac r{2^j}\right)
\]
for any \(r\geqslant 0\).
Splitting the sum on the right into two ranges, \(j \leqslant \log_2 r\) and \(j > \log_2 r\), we apply the condition of the theorem~\eqref{eq:functional-equation-E-condition} to the first range, and use the fact that \(E\bigl(g^{(N+1)}, r\bigr)\) remains bounded for bounded \(r\) in the second range. This yields
\[
\begin{split}
E\left(f^{\left( N+1\right) };r\right)&\leqslant\sum_{0\leqslant j \leqslant \log_2 r}\frac{1}{2^{jN}}E\left(g^{\left( N+1\right) };\frac r{2^j}\right)+\sum_{j>\log_2 r}\frac{1}{2^{jN}}E\left(g^{\left( N+1\right) };\frac r{2^j}\right)
\\
&=O\left(\sum_{0\leqslant j \leqslant \log_2 r}\frac{1}{2^{jN}}\left(\frac{1}{{r}/{2^j}}\right)^N+\sum_{j>\log_2 r}\frac{1}{2^{jN}}\right).
\end{split}
\]
Thus, for \(N\geqslant1\),
	\[
	E(f^{(N+1)};r)=O\left(\frac{\log r}{r^N}\right)
	\]
	as \(r\to +\infty\). Together with the \(N=0\) estimate, this proves the same bound for every \(N\geqslant0\). Applying this estimate to evaluate the right side of the inequality of Theorem \ref{thm:main-depoissonization} completes the proof of the theorem.
\end{proof}

\section{Applications}
\label{sec:applications}
The first three applications begin with first-order depoissonization
estimates. The Stirling-number application also illustrates the
higher-order bound of Theorem~\ref{thm:main-depoissonization}, while the final
symmetric-trie application combines the monotonicity estimate with the
functional-equation bound from Theorem~\ref{thm:functional-equation-bound}.

\subsection{Summatory sequences}
The first-order estimate in Theorem~\ref{thm:first-order-bound} has the following consequence for summatory sequences.
\begin{cor} Suppose  \(a_0,a_1,\ldots\) is a sequence of nonnegative numbers such that the series defining the function
\[
g(t)=e^{-t}\sum ^{\infty }_{n=0}\frac{a_{n}}{n!}t^{n}
\]
has an infinite radius of convergence. Then
\[
\left|a_0+a_1+\cdots+a_{n-1}-\int_0^{n}g(t)\,dt\right|\leqslant 2\sqrt{n}g(n)
\]
for all \(n\geqslant 1\).
\end{cor}
\begin{proof}
Set \(A_0=0\), \(A_n=a_0+\cdots+a_{n-1}\) for \(n\geqslant1\),
and consider the exponential generating series
\[
S(x)=\sum_{n\geqslant0}A_n\frac{x^n}{n!}.
\]
Since \(A_{n+1}-A_n=a_n\), coefficientwise differentiation gives
\[
S'(x)-S(x)
=\sum_{n\geqslant0}a_n\frac{x^n}{n!}
=e^xg(x).
\]
The right-hand side is entire, and the solution with \(S(0)=0\) is
entire. Thus the series defining \(S\) has infinite radius of
convergence, and
\[
\bigl(e^{-x}S(x)\bigr)'=g(x),
\qquad
e^{-x}S(x)=\int_0^x g(t)\,dt.
\]
Hence the sequence \((A_n)\) has Poisson transform
\(f(x)=\int_0^xg(t)\,dt\). Since
\(\Delta A_m=a_m\geqslant0\), inequality~\eqref{eq:first-order-bound-diagonal} applies and
gives the result, with \(E(f';r)=f'(r)=g(r)\).

\end{proof}

\subsection{Poissonized and fixed-index expectations}
Proposition~\ref{prop:poissonized-expectation-bound} applies Corollary~\ref{cor:monotone-depoissonization} to the sequence \(A_m=\mathbb{E}X_m\). The derivative of the corresponding Poisson transform is expressed as a covariance and controlled by the variance of the poissonized variable.

Suppose
\[
X_0, X_1, X_2, \ldots
\]
is a sequence of random variables, and let \(\eta(\lambda)\) be a Poisson random variable with parameter \(\lambda > 0\), i.e.,
\[
\mathbb{P}\bigl(\eta(\lambda) = m\bigr) = e^{-\lambda}\frac{\lambda^m}{m!},
\qquad m = 0,1,2,\ldots
\]
Assume further that \(\eta(\lambda)\) is independent of the sequence \(X_0, X_1, X_2, \ldots\).

In many applications (see, e.g., \citet{hwang_janson_2008}), the \emph{poissonized} random variable
\[
X_{\eta(n)}
\]
is more tractable than the fixed-index variable \(X_n\): for example, it may admit explicit moment formulas or have a simpler limiting distribution. Since \(\eta(n)\) has mean \(n\) and standard deviation \(\sqrt n\),
the random index \(\eta(n)\) fluctuates around \(n\) on the scale
\(\sqrt n\). This suggests that, whenever the distributions of the
random variables \(X_j\) vary sufficiently slowly with \(j\) on this
scale, the poissonized variable and the fixed-index variable may be
close in distribution. In particular, one may expect
\[
X_{\eta(n)}\approx X_n.
\]
The following proposition quantifies the proximity between the expectations of the poissonized and unpoissonized variables, under the assumption that the expectation sequence is monotone. The bound is expressed in terms of the variance of the poissonized variable.
\begin{prop}
\label{prop:poissonized-expectation-bound}
Assume that \(X_{\eta(R)}\in L^2\) for every \(R>0\), and suppose the sequence of expectations
\[
\mathbb{E}X_0, \; \mathbb{E}X_1, \; \mathbb{E}X_2, \;\ldots
\]
is monotone (either nondecreasing or nonincreasing). Then, for every integer \(n\geqslant1\) and every \(R > 0\),
\[
\bigl| \,\mathbb{E}X_n - \mathbb{E}X_{\eta(R)} \,\bigr|
 \;\leqslant\;
 \frac{e^R}{R^n}\,\frac{n!}{\sqrt{R}} \,
 \sqrt{ \mathrm{Var}\!\left( X_{\eta(R)}\right) }.
\]
In particular, when \(R=n\) one obtains the sharper bound
\[
\bigl| \,\mathbb{E}X_n - \mathbb{E}X_{\eta(n)} \,\bigr|
 \;\leqslant\;
 2 \sqrt{ \mathrm{Var}\!\left( X_{\eta(n)}\right) }.
\]
\end{prop}

\begin{proof}[Proof of Proposition \ref{prop:poissonized-expectation-bound}]
	The \(L^2\) assumption ensures that all expectations, variances, and covariances below are finite.
	We will apply Corollary \ref{cor:monotone-depoissonization} with \[
A_m=\mathbb{E}X_m.
\]
For every \(r>0\), put \(p_m(r)=e^{-r}r^m/m!\). By independence,
the \(L^2\) assumption, the Cauchy--Schwarz inequality, and Jensen's
inequality,
\[
\begin{aligned}
e^{-r}\sum_{m\geqslant0}|\mathbb{E}X_m|\frac{r^m}{m!}
&=\sum_{m\geqslant0}p_m(r)|\mathbb{E}X_m|\\
&\leqslant
\left(\sum_{m\geqslant0}p_m(r)(\mathbb{E}X_m)^2\right)^{1/2}\\
&\leqslant
\left(\sum_{m\geqslant0}p_m(r)\mathbb{E}X_m^2\right)^{1/2}\\
&=\|X_{\eta(r)}\|_2<\infty.
\end{aligned}
\]
Thus \(\sum_{m\geqslant0}\mathbb{E}X_m z^m/m!\) converges absolutely
for every \(z\in\mathbb C\), and the corresponding Poisson transform
is entire.
Then the generating function
\[
f(\lambda)=e^{-\lambda}\sum_{m=0}^\infty\mathbb{E}X_m\frac{\lambda^m}{m!}
\]
can be viewed as an expectation
\[
f(\lambda)=\mathbb{E}X_{\eta(\lambda)}
\]
and so its derivative can be expressed in terms of expectations involving the Poisson random variable \(\eta(\lambda)\) as
\[
\begin{aligned}
f'\left( \lambda \right) &=\frac{1}{\lambda}\sum ^{\infty }_{m=0}m\mathbb{E}X_{m}e^{-\lambda }\frac{\lambda ^{m}}{m!}-\sum ^{\infty }_{m=0}\mathbb{E}X_{m}e^{-\lambda }\frac{\lambda ^{m}}{m!}
=\frac{1}{\lambda }\mathbb{E}\eta \left( \lambda \right) X_{\eta\left( \lambda \right) }-\mathbb{E}X_{\eta \left( \lambda\right)}.
\end{aligned}
\]
We can further express the derivative \(f'\left( \lambda \right)\) as a covariance
\[
\begin{aligned}f'\left( \lambda\right) =\frac{1}{\lambda }\mathbb{E}\left(\eta \left( \lambda \right) -\lambda \right)X_{\eta\left( \lambda \right) }
=\frac{1}{\lambda }\mathbb{E}(\eta \left( \lambda \right) -\lambda )\left(X_{\eta\left( \lambda \right) } -\mathbb{E}X_{\eta(\lambda)}\right)=\frac1\lambda\mathrm{Cov}\bigl(\eta \left( \lambda \right), X_{\eta(\lambda)}\bigr).
\end{aligned}
\]
The Cauchy--Schwarz inequality now gives
\[
|f'(\lambda)|\leqslant \frac{1}{\lambda}\sqrt{\mathrm{Var}\bigl(\eta(\lambda)\bigr)}\sqrt{\mathrm{Var}\bigl(X_{\eta(\lambda)}\bigr)}=\frac{\sqrt{\mathrm{Var}\bigl(X_{\eta(\lambda)}\bigr)}}{\sqrt{\lambda}}
\]
since \(\mathrm{Var}\bigl(\eta(\lambda)\bigr)=\lambda\).
Substituting this bound for \(f'(\lambda)\), together with the expressions for \(A_n\) and \(f(\lambda)\), into Corollary \ref{cor:monotone-depoissonization} completes the proof.
\end{proof}

\subsection{Stirling numbers of the second kind}

The Stirling number of the second kind
\[
\Stirling{m}{k}
\]
is the number of partitions of an \(m\)-element set into \(k\)
nonempty blocks. Its exponential generating function is
\[
\sum_{m=0}^{\infty}\Stirling{m}{k}\frac{z^m}{m!}
=
\frac{(e^z-1)^k}{k!},
\qquad k\geqslant 1;
\]
see, for example, \citet[p.~206]{comtet_1974}.

The following proposition gives a uniform nonasymptotic approximation
as a direct application of Corollary~\ref{cor:monotone-depoissonization}.

\begin{prop}
\label{prop:stirling-first-order}
For all integers \(n\geqslant k\geqslant 1\),
\begin{equation}
\label{eq:stirling-absolute-bound}
\left|
\Stirling{n}{k}\frac{k!}{k^n}
-
\left(1-e^{-n/k}\right)^k
\right|
\leqslant
2\sqrt{n}\,e^{-n/k}
\left(1-e^{-n/k}\right)^{k-1}.
\end{equation}
Equivalently,
\begin{equation}
\label{eq:stirling-relative-bound}
\left|
\frac{k!\Stirling{n}{k}}
{k^n\left(1-e^{-n/k}\right)^k}
-1
\right|
\leqslant
\frac{2\sqrt{n}\,e^{-n/k}}
{1-e^{-n/k}}.
\end{equation}
\end{prop}

\begin{proof}
For a fixed integer \(k\geqslant 1\), define
\[
A_m^{(k)}
=
\Stirling{m}{k}\frac{k!}{k^m},
\qquad m\geqslant 0.
\]
The quantity \(A_m^{(k)}\) is the probability that, after \(m\)
independent balls have been placed uniformly into \(k\) labelled
boxes, no box is empty. Under the natural coupling obtained by adding
one further ball, this event is increasing in \(m\). Hence the sequence
\(\bigl(A_m^{(k)}\bigr)_{m\geqslant0}\) is nondecreasing.

The same monotonicity can also be verified algebraically.
Differentiating the exponential generating function above and
comparing coefficients gives the standard recurrence
\[
\Stirling{m+1}{k}
=
k\Stirling{m}{k}+\Stirling{m}{k-1}.
\]
Consequently, for \(k\geqslant2\),
\[
A_{m+1}^{(k)}-A_m^{(k)}
=
\left(1-\frac{1}{k}\right)^m A_m^{(k-1)}
\geqslant0.
\]
For \(k=1\), one has \(A_0^{(1)}=0\) and
\(A_m^{(1)}=1\) for every \(m\geqslant1\), so monotonicity also holds
in this case.

Using the exponential generating function of the Stirling numbers, the
Poisson transform of this sequence is
\[
\begin{aligned}
f_k(x)
&=
e^{-x}\sum_{m=0}^{\infty}
A_m^{(k)}\frac{x^m}{m!} \\
&=
e^{-x}\sum_{m=0}^{\infty}
\Stirling{m}{k}\frac{k!}{m!}
\left(\frac{x}{k}\right)^m \\
&=
e^{-x}\left(e^{x/k}-1\right)^k
=
\left(1-e^{-x/k}\right)^k.
\end{aligned}
\]
Consequently,
\[
f_k'(x)
=
e^{-x/k}\left(1-e^{-x/k}\right)^{k-1}.
\]
Applying inequality~\eqref{eq:monotone-diagonal-bound} of
Corollary~\ref{cor:monotone-depoissonization} to \(A_m^{(k)}\), with \(R=n\),
gives
\[
\left|A_n^{(k)}-f_k(n)\right|
\leqslant
2\sqrt{n}\,f_k'(n),
\]
which is precisely~\eqref{eq:stirling-absolute-bound}. Dividing by
\(\left(1-e^{-n/k}\right)^k\) gives
\eqref{eq:stirling-relative-bound}.
\end{proof}
As an immediate consequence, for every fixed \(\varepsilon>0\),
\begin{equation}
\label{eq:stirling-asymptotic-range}
\Stirling{n}{k}
=
\frac{k^n}{k!}
\left(1-e^{-n/k}\right)^k
\left(1+O_{\varepsilon}\!\left(n^{-\varepsilon/2}\right)\right)
\end{equation}
uniformly, as \(n\to\infty\), for integers \(k\) satisfying
\[
1\leqslant k
\leqslant
\frac{2n}{(1+\varepsilon)\log n}.
\]
Indeed, in this range
\[
e^{-n/k}
\leqslant
n^{-(1+\varepsilon)/2},
\]
so the right-hand side of
\eqref{eq:stirling-relative-bound} is
\(O_{\varepsilon}(n^{-\varepsilon/2})\).

For a broader asymptotic analysis of the Stirling numbers of the second
kind, including higher-order expansions in the polynomials \(\tau_j\)
obtained by an elementary finite-differencing analysis of the
inclusion--exclusion formula, see
\citet{hwang_li_zacharovas_2026}. The present approach gives an alternative derivation. Proposition~\ref{prop:stirling-first-order}
uses monotonicity to obtain the leading approximation
\((1-e^{-n/k})^k\). This can be further extended using Theorem~\ref{thm:main-depoissonization}, together with
the properties of \(E\) to yield higher-order expansions. We illustrate
this with the  expansion corresponding to \(N=2\).

\begin{prop}[Second-order expansion for Stirling numbers]
\label{prop:stirling-second-order}
For every integer \(k\geqslant3\), let
\[
A_n^{(k)}:=\frac{k!}{k^n}\Stirling{n}{k},
\qquad
u_{n,k}:=\frac{\sqrt n\,e^{-n/k}}{1-e^{-n/k}}.
\]
Then, for every integer \(n\geqslant2\),
\[
\left|
\frac{A_n^{(k)}}{(1-e^{-n/k})^k}
-1
-\frac{n}{2k}
\frac{e^{-n/k}(1-ke^{-n/k})}{(1-e^{-n/k})^2}
\right|
\leqslant
33
\left(
\frac{n}{k^2}u_{n,k}
+3\frac{\sqrt n}{k}u_{n,k}^{\,2}
+u_{n,k}^{\,3}
\right).
\]
\end{prop}

\begin{proof}
Let
\[
f_k(x):=(1-e^{-x/k})^k.
\]
As in the proof of Proposition~\ref{prop:stirling-first-order},
\(f_k\) is the Poisson transform of
\((A_m^{(k)})_{m\geqslant0}\). Since
\[
\tau_0(n)=1,\qquad
\tau_1(n)=0,\qquad
\tau_2(n)=-n,
\]
Theorem~\ref{thm:main-depoissonization} with \(N=2\) gives
\begin{equation}
\label{eq:stirling-second-order-depoissonization-bound}
\left|
A_n^{(k)}
-f_k(n)+\frac n2 f_k''(n)
\right|
\leqslant
33n^{3/2}E(f_k^{(3)};n).
\end{equation}

We next estimate the \(E\)-term on the right. Set
\[
h_p(x):=
e^{-px/k}(1-e^{-x/k})^{k-p},
\qquad 1\leqslant p\leqslant3.
\]
For \(p<k\), differentiation gives
\[
h_p'(x)
=
-\frac pk h_p(x)+\frac{k-p}{k}h_{p+1}(x).
\]
Since \(f_k'(x)=h_1(x)\), two applications of this identity yield
\[
\begin{aligned}
f_k''(x)
&=
-\frac1k h_1(x)+\frac{k-1}{k}h_2(x),\\
f_k^{(3)}(x)
&=
\frac1{k^2}h_1(x)
-\frac{3(k-1)}{k^2}h_2(x)
+\frac{(k-1)(k-2)}{k^2}h_3(x).
\end{aligned}
\]

For \(1\leqslant p<k\), we have
\[
\begin{aligned}
h_p(x)
&=
e^{-x}(e^{x/k}-1)^{k-p}\\
&=
e^{-x}\sum_{m=0}^{\infty}
\frac{(k-p)!}{k^m}
\Stirling{m}{k-p}\frac{x^m}{m!}.
\end{aligned}
\]
When \(k=p=3\), the corresponding endpoint is
\(h_3(x)=e^{-x}\), whose Poisson coefficients are also nonnegative.
Thus every function \(h_p\), \(1\leqslant p\leqslant3\), occurring
above has nonnegative Poisson coefficients, and consequently
\[
E(h_p;n)=h_p(n).
\]
The linear-combination inequality
\eqref{eq:E-linear-combination} therefore gives
\begin{equation}
\label{eq:stirling-third-derivative-E-bound}
\begin{aligned}
E(f_k^{(3)};n)
\leqslant{}&
\frac{e^{-n/k}(1-e^{-n/k})^{k-1}}{k^2}
\\
&+\frac{3(k-1)}{k^2}
e^{-2n/k}(1-e^{-n/k})^{k-2}
\\
&+\frac{(k-1)(k-2)}{k^2}
e^{-3n/k}(1-e^{-n/k})^{k-3}.
\end{aligned}
\end{equation}

It remains to combine
\eqref{eq:stirling-second-order-depoissonization-bound} and
\eqref{eq:stirling-third-derivative-E-bound}. Put
\[
q:=e^{-n/k}.
\]
Then \(0<q<1\), and
\[
f_k(n)=(1-q)^k,
\qquad
f_k''(n)
=
-\frac1k q(1-q)^{k-2}(1-kq).
\]
Dividing
\eqref{eq:stirling-second-order-depoissonization-bound} by
\(f_k(n)>0\), substituting
\eqref{eq:stirling-third-derivative-E-bound}, and then simplifying
gives
\[
\begin{aligned}
&\left|
\frac{A_n^{(k)}}{(1-q)^k}
-1
-\frac{n}{2k}\frac{q(1-kq)}{(1-q)^2}
\right|
\\
&\quad\leqslant
\frac{33n^{3/2}}{(1-q)^k}E(f_k^{(3)};n)
\\
&\quad\leqslant
33n^{3/2}
\left(
\frac{q}{k^2(1-q)}
+\frac{3(k-1)q^2}{k^2(1-q)^2}
+\frac{(k-1)(k-2)q^3}{k^2(1-q)^3}
\right)
\\
&\quad\leqslant
33n^{3/2}
\left(
\frac{q}{k^2(1-q)}
+\frac{3q^2}{k(1-q)^2}
+\frac{q^3}{(1-q)^3}
\right)
\\
&\quad=
33
\left(
\frac{n}{k^2}u_{n,k}
+3\frac{\sqrt n}{k}u_{n,k}^{\,2}
+u_{n,k}^{\,3}
\right),
\end{aligned}
\]
where the last inequality uses
\[
\frac{3(k-1)}{k^2}\leqslant\frac3k,
\qquad
\frac{(k-1)(k-2)}{k^2}\leqslant1,
\]
and the last equality follows from
\[
u_{n,k}=\frac{\sqrt n\,q}{1-q}.
\]
This is the asserted inequality.
\end{proof}

\subsection{Symmetric tries}

\begin{exm}[Application to symmetric tries]
A straightforward example of an application of Theorem~\ref{thm:functional-equation-bound} concerns the expectation of the size \(S_n\) of a random symmetric trie. The depoissonization step in this problem is usually treated using complex analytic tools as in \citet{jacquet_szpankowski_1998}. Here we present a way to obtain the depoissonization expansion relying only on elementary tools developed in this paper.

\textbf{First-term expansion (using monotonicity of \(\mathbb{E}S_m\)).}
The generating function of  expectations  of \(S_n\) satisfies the functional equation
\[
h(x) = 2h\!\left(\tfrac{x}{2}\right) + e^{-x}(e^x - 1 - x),
\]
where
\[
h(x) = e^{-x}\sum_{n=0}^\infty \frac{\mathbb{E}S_n}{n!}\,x^n, 
\qquad \mathbb{E}S_0 = \mathbb{E}S_1 = 0.
\]

Differentiating the functional equation gives
\[
h'(x) = h'\!\left(\tfrac{x}{2}\right) + x e^{-x}.
\]
Equivalently,
\[
e^{-x}\sum_{n=0}^\infty \frac{\Delta \mathbb{E}S_n}{n!}\,x^n
= e^{-x/2}\sum_{n=0}^\infty \frac{\Delta \mathbb{E}S_n}{n!}\,\Bigl(\frac{x}{2}\Bigr)^n
+ x e^{-x}.
\]
Multiplying through by \(e^x\) yields
\[
\sum_{n=0}^\infty \frac{\Delta \mathbb{E}S_n}{n!}\,x^n
= \sum_{m=0}^\infty \frac{1}{m!}\Bigl(\frac{x}{2}\Bigr)^m 
   \sum_{n=0}^\infty \frac{\Delta \mathbb{E}S_n}{n!}\,\Bigl(\frac{x}{2}\Bigr)^n 
   + x.
\]
 Comparing coefficients of \(x^n\) on both sides gives
\[
\Delta \mathbb{E}S_n
= \frac{1}{2^n}\sum_{j=0}^n \binom{n}{j}\,\Delta \mathbb{E}S_j + \delta_{n,1},
\]
where \(\delta_{n,1}=1\) if \(n=1\) and \(\delta_{n,1}=0\) otherwise.  
Equivalently,
\[
\Delta \mathbb{E}S_n
= \frac{1}{2^n-1}\sum_{j=0}^{n-1} \binom{n}{j}\,\Delta \mathbb{E}S_j + 2\delta_{n,1}, 
\qquad n \geqslant 1.
\]
Since \(\delta_{n,1}\geqslant 0\) and \(\Delta \mathbb{E}S_0 
= \mathbb{E}S_1 - \mathbb{E}S_0 = 0\), one easily shows by induction that
\[
\Delta \mathbb{E}S_n \geqslant 0, \qquad n \geqslant 0.
\]
Therefore Corollary~\ref{cor:monotone-depoissonization} applies and yields
\[
\bigl| \mathbb{E}S_n - h(n) \bigr| \leqslant 2\sqrt{n}\,\bigl|h'(n)\bigr|.
\]

\textbf{Full asymptotic expansion (using properties of the operator \(E\)).}
Extending this argument to higher-order expansions would require \(\Delta^{N+1}\mathbb{E}S_n\) to have constant sign for some \(N>0\), which is not the case (at least for small values of \(n\)). Instead, we apply Theorem~\ref{thm:functional-equation-bound} to \(f=h\). The
corresponding function \(g\) is
\[
g(z):=h(z)-2h\left(\frac z2\right)
     =e^{-z}(e^z-1-z).
\]
The initial conditions give \(h(0)=h'(0)=0\), and iteration of the
functional equation, coefficientwise at the origin, gives
\[
h(z)=\sum_{j\geqslant0}2^j g\left(\frac{z}{2^j}\right).
\]
This representation also proves the required entireness. Indeed,
\(g(z)=O(z^2)\) as \(z\to0\). Hence, on every compact set and for all
sufficiently large \(j\),
\[
2^j g\left(\frac{z}{2^j}\right)=O(2^{-j})
\]
uniformly in \(z\). The series therefore converges locally uniformly
and defines an entire function. It satisfies the functional equation
and the normalization \(h(0)=h'(0)=0\), which determine its power-series
coefficients, so it is the Poisson transform \(h\) above. In particular,
\(h\) is entire, as required in Theorem~\ref{thm:functional-equation-bound}.
In this setting,
\[
g^{(N+1)}(x) = (-1)^{N+1} e^{-x}(N-x),
\]
so that
\[
E\bigl(g^{(N+1)};r\bigr) = e^{-r}(N+r), \qquad r>0.
\]
Since this quantity decays exponentially, condition~\eqref{eq:functional-equation-E-condition} is satisfied for every \(N \geqslant 0\). Consequently,
\[
\mathbb{E}S_n - \sum_{m=0}^N \frac{h^{(m)}(n)}{m!}\,\tau_m(n)
= O\!\left(\frac{\log n}{n^{\frac{N-1}{2}}}\right),
\qquad n \to +\infty,
\]
for any fixed \(N \geqslant 0\).
\end{exm}

\section*{Acknowledgements}

The author wrote this paper while serving as a visiting associate professor at the Institute of Statistical Science, Academia Sinica (Taiwan). He sincerely thanks Prof. Hsien-Kuei Hwang for his generous hospitality during the visit and also for drawing the author's attention to Ramanujan's expansions as presented in Entry 10 of  \citet{berndt_1985} and to the paper by \citet{yu_2009}. The author is also grateful to the anonymous referee for a careful reading of the manuscript and for numerous constructive comments and suggestions. These led to substantial improvements in the organization and exposition of the paper and to clearer formulations of several assumptions and results.

\par\medskip
\noindent\textit{Use of generative AI.}
During the revision of this manuscript, the author used OpenAI's
generative AI tools  to assist with manuscript
organization and exposition and with checking mathematical arguments
and consistency. All suggestions and edits were critically reviewed
and verified by the author, who takes full responsibility for the
final manuscript.

\begin{appendices}
\renewcommand*{\theHequation}{appendix.\arabic{section}.\arabic{equation}}
\section{Properties of polynomials \texorpdfstring{$\tau_m(n)$ and $C_m(\lambda,n)$}{tau\_m(n) and C\_m(lambda,n)}}

\label{app:polynomial-properties}

Let us recall that, initially for \(n\in\mathbb N_0\), we have defined
\(\tau_m(n)\) by
\[
\tau_m(n) = n! [z^n] ( z-n)^m e^z.
\]

\begin{prop}
For every real \(n\), if \(|x|<1\), then
\begin{equation}
	\label{eq:tau-generating-function-appendix}
	\sum ^{\infty }_{m=0}\frac{\tau_{m}\left( n\right) }{m!}x^m=\bigl((1+x)e^{-x}\bigr)^n
\end{equation}	
and as a consequence
\begin{equation}
\label{eq:tau-derivative-formula}
\tau_{m}\left( n\right) =\frac{d^m}{dx^m}\bigl((1+x)e^{-x}\bigr)^n\Bigg|_{x=0}.
\end{equation}
\end{prop}
\begin{proof}
First let \(n\in\mathbb N_0\). To derive the generating function, we apply the expansion
\eqref{eq:charlier-poisson-expansion} to the function $f(z) = e^{zx}$.
In this case,
\[
f(z) = e^{zx} = e^{-z} e^{z(1+x)} 
       = e^{-z}\sum_{n=0}^{\infty} \frac{(1+x)^n}{n!}\,z^n.
\]
Thus the coefficients are $A_n = (1+x)^n$, and substituting into 
\eqref{eq:charlier-poisson-expansion} gives
\[
(1+x)^n = \sum_{m=0}^{\infty} \frac{e^{nx}x^m}{m!}\,\tau_m(n).
\]
Multiplying both sides by $e^{-nx}$ yields
\eqref{eq:tau-generating-function-appendix}. The identity \eqref{eq:tau-derivative-formula} follows
by differentiating identity~\eqref{eq:tau-generating-function-appendix} \(m\) times at
$x=0$. For each fixed \(m\), the derivative on the right-hand side of
\eqref{eq:tau-derivative-formula} is a polynomial in \(n\). It therefore supplies
the polynomial continuation of \(\tau_m(n)\) from \(n\in\mathbb N_0\)
to \(n\in\mathbb R\), for which the same local generating-function
identity holds.
\end{proof}

Thus the formulas in the preceding proposition show that
\(\tau_m(n)\) is a polynomial in \(n\) and extend its coefficient
definition from nonnegative integer \(n\) to every real \(n\).
\begin{prop}
\label{prop:tau-bound}
	Functions \((-1)^m\tau_m(-y)\) are polynomials in the variable \(y\) with nonnegative coefficients of degree at most \(\lfloor{m/2}\rfloor\), and the inequality
	\[
	|\tau_m(y)|\leqslant m!\frac{|y|^{\lfloor{m/2}\rfloor}}{2}
	\]
	holds for all \(m\geqslant 1\) and \(|y|\geqslant 1\).
\end{prop}
\begin{proof} Replacing \(x\) and \(n\) with \(-x\) and \(-y\), respectively, in identity~\eqref{eq:tau-generating-function-appendix} and expanding it into a Taylor series with respect to \(y\), we get
\begin{equation}
	\sum ^{\infty }_{m=0}\frac{(-1)^m\tau_{m}\left( -y\right) }{m!}x^m=\big((1-x)e^x\big)^{-y}
=\exp\!\Big(y\sum_{r\geqslant2}\frac{x^r}{r}\Big)=\sum_{k=0}^\infty \frac{y^k}{k!}\left(\sum_{r=2}^\infty \frac{x^r}{r}\right)^k.
\end{equation}
This implies that \((-1)^m\tau_m(-y)\) is a polynomial with nonnegative coefficients in the variable \(y\) of degree at most \(\lfloor{m/2}\rfloor\). This means that
to evaluate \(\tau_m(y)\) for \(y\geqslant 0\) we can use the upper bound
\[
|\tau_m(y)|\leqslant (-1)^m\tau_m(-y).
\]
Since all the coefficients of the polynomial \((-1)^m\tau_m(-y)\) are nonnegative, for \(y\geqslant1\) its value can be bounded by \(y^{\lfloor m/2\rfloor}\) multiplied by the sum of the coefficients, which equals \((-1)^m\tau_m(-1)\). Thus we have an inequality
\begin{equation}
\label{eq:tau-negative-bound}
0\leqslant(-1)^m\tau_m(-y)\leqslant 	(-1)^m\tau_m(-1)y^{\lfloor{m/2}\rfloor}
\end{equation}
that holds for all \(y\geqslant 1\). To compute the value of \(\tau_m(-1)\) we insert \(n=-1\)  into the identity (\ref{eq:tau-generating-function-appendix}) and after replacing \(x\) with \(-x\) we obtain
\[
\sum ^{\infty }_{m=0}\frac{(-1)^m\tau_{m}\left( -1\right) }{m!}x^m=\frac{1}{(1-x)e^x
}=\sum_{k=0}^\infty x^k\sum_{j=0}^\infty \frac{(-x)^j}{j!}=\sum_{m=0}^\infty \left(\sum_{j\leqslant m}\frac{(-1)^j}{j!}\right)x^m
\]
hence
\[
\frac{(-1)^m\tau_{m}\left( -1\right)}{m!} =\sum_{j\leqslant m}\frac{(-1)^j}{j!}.
\]
Inserting this  expression into the inequality (\ref{eq:tau-negative-bound}) and noting that the sum on the right side of the above identity does not exceed \(1/2\) for all \(m\geqslant 1\) we obtain the proof of the inequality of our proposition.
\end{proof}
An alternative derivation of a slightly weaker estimate, 
\(|\tau_m(n)| \leqslant m!n^{m/2}\) for integers \(m,n \geqslant 0\), 
was previously obtained  by the author and the coauthors in~\citet{hwang_li_zacharovas_2026},
using a different method.

The following lemma states an interesting identity that is an immediate
consequence of Parseval's identity for Charlier polynomials. We give an
independent proof based only on the properties of the polynomials
\(\tau_m\).
\begin{lem} For every positive integer \(n\),
\[
\sum_{m=0}^\infty\frac{\tau_m^2(n)}{m!n^m}=\frac{n!}{\left(\frac{n}{e}\right)^n}.
\]	
\end{lem}
\begin{proof} Recall that \eqref{eq:tau-coefficient-definition} defines \(\tau_m(n)\) as the coefficient of \(z^n\) in \(n!(z-n)^m e^z\). Therefore, we can represent \(\tau_m^2(n)\) as a product of the coefficients of \(z^n\) and \(w^n\) in \(n!(z-n)^m e^z\) and \(n!(w-n)^m e^w\), respectively. Substituting this representation into the sum on the left-hand side gives
\[
\begin{split}
\sum_{m=0}^\infty\frac{\tau_m^2(n)}{m!n^m}&=\sum_{m=0}^\infty\frac{n! [z^n] (z-n)^m e^zn! [w^n] (w-n)^m e^w}{m!n^m}
\\
&=n!^2 [z^n][w^n] e^{z+w}\sum_{m=0}^\infty\frac{(z-n)^m   (w-n)^m }{m!n^m}
\\
&=n!^2 [z^n][w^n] e^{z+w +\frac{(z-n)(w-n)}{n}}
\\
&=n!^2 [z^n][w^n] e^{n +\frac{zw}{n}}.
\end{split}
\]
Since
\[
[z^n][w^n]\,e^{zw/n}
= \frac{1}{n!}\left(\frac{1}{n}\right)^{n}
\]
we finally obtain the identity of the lemma.
\end{proof}

By bounding an individual summand in the identity of the preceding
lemma by the value of the full sum, we obtain the following estimate,
which can be better than the estimate of Proposition \ref{prop:tau-bound}
if \(m\) is large enough.
\begin{cor}
For all integers \(n\geqslant1\) and \(m\geqslant0\),
\[
{|\tau_m(n)|}\leqslant \sqrt{\frac{n!}{\left(\frac{n}{e}\right)^n}}\sqrt{m!}{n^{m/2}}
\]
\end{cor}

For \(\lambda\ne0\), our definition of the Charlier polynomials \(C_m(\lambda,n)\) in \eqref{eq:charlier-coefficient-definition} means that we can alternatively define them by the identity\begin{equation}
\label{eq:charlier-coefficient-generating-function}
\sum _{n\geqslant 0}\frac{C_{m}\left( \lambda ,n\right)}{n!}z^{n}=\left( \frac{z}{\lambda }-1\right) ^{m}e^{z}.	
\end{equation}
Note that this coefficient definition is initially formulated only for
nonnegative integer values of \(n\). To extend the definition to real
\(n\), we use the following classical generating function. This
generating function and the finite-sum formula in the corollary below
are standard; see, for example,
\citet[Sec.~2.81, Eq.~(2.81.2), pp.~34--35]{szego_1975},
\citet[pp.~170--172]{chihara_1978}, and
\citet[Sec.~9.14, Eq.~(9.14.1)]{koekoek_lesky_swarttouw_2010}.
We include short proofs to fix our normalization and to keep the
appendix self-contained.
\begin{thm}
\label{thm:charlier-generating-function}
	For every real \(\lambda\ne0\) and every real \(n\),
	\begin{equation}
	\label{eq:charlier-generating-function}
	\sum _{m\geqslant 0}\frac{w^m}{m!} C_{m}\left( \lambda ,n\right)
= e^{-w} \left( \frac{w}{\lambda }+1\right)^n,
\qquad |w|<|\lambda|,
	\end{equation}
where the branch is determined by the Taylor expansion at \(w=0\).
If \(n\in\mathbb N_0\), the right-hand side is entire in \(w\), and
the identity extends to every \(w\in\mathbb C\).
\end{thm}
\begin{proof}
First let \(n\in\mathbb N_0\). Multiply both sides of the generating
function~\eqref{eq:charlier-coefficient-generating-function} by \(w^m/m!\) and sum over
\(m\) from zero to infinity. As an identity of power series, this gives
\[
\sum _{m\geqslant 0}\frac{w^m}{m!}
\sum _{r\geqslant 0}\frac{C_{m}\left( \lambda ,r\right)}{r!}z^{r}= \sum _{m\geqslant 0}\frac{w^m}{m!}\left( \frac{z}{\lambda }-1\right) ^{m}e^{z}= e^{z}e^{w\left( \frac{z}{\lambda }-1\right) }= e^{-w}e^{z\left( \frac{w}{\lambda }+1\right) }
\]
Taking the coefficient of \(z^n/n!\) proves the identity for
\(n\in\mathbb N_0\), for all \(w\). For fixed \(m\), the coefficient
of \(w^m/m!\) in the Taylor expansion at zero of
\(e^{-w}(1+w/\lambda)^n\) is a polynomial in \(n\). These coefficients
agree with \(C_m(\lambda,n)\) at every nonnegative integer \(n\), and
hence provide its polynomial continuation to real \(n\). For real
\(n\), the Taylor series of the right-hand side converges for
\(|w|<|\lambda|\), with this branch and these continued coefficients.
This proves~\eqref{eq:charlier-generating-function} in the stated domain.
\end{proof}
The identity of the above theorem can be used to obtain an explicit expression for Charlier polynomials, formulated in the next corollary, that
is valid not only for integer values of \(n\) but also for real ones.
\begin{cor}
For \(\lambda\ne0\), \(m\in\mathbb N_0\), and \(x\in\mathbb R\),
\[
C_m(\lambda,x)=\sum_{j=0}^m\left(-1\right)^{m-j}\binom{m}{j}\frac{x(x-1)\cdots(x-j+1)}{\lambda^j}.
\]	
\end{cor}
\begin{proof}
For \(|w|<|\lambda|\), Theorem~\ref{thm:charlier-generating-function} gives
\[
\sum _{m\geqslant 0}\frac{w^m}{m!} C_{m}\left( \lambda ,x\right)
= e^{-w} \left( \frac{w}{\lambda }+1\right)^x=\sum_{m=0}^\infty\frac{(-w)^m}{m!}
\sum_{j=0}^\infty\frac{x(x-1)\cdots(x-j+1)}{j!}\left(\frac{w}{\lambda}\right)^j
\]
Comparing the coefficients of \(w^m\) yields the explicit expression
stated in the corollary.
\end{proof}

The formula above agrees with the classical expression in
\citet[Sec.~2.81, Eq.~(2.81.2)]{szego_1975}, after accounting for
normalization. If \(\lambda>0\), Szeg\H{o}'s orthonormal
Poisson--Charlier polynomial \(p_m(x)\) with parameter \(\lambda\)
satisfies
\[
p_m(x)
=
\left(\frac{\lambda^m}{m!}\right)^{1/2}
C_m(\lambda,x).
\]

To compare our convention with that of
\citet[\S9.14]{koekoek_lesky_swarttouw_2010}, let
\(C_m^{\mathrm{KLS}}(x;\lambda)\) denote the Charlier polynomial in
their normalization. By their hypergeometric representation
\citep[p.~247, (9.14.1)]{koekoek_lesky_swarttouw_2010},
\[
C_m^{\mathrm{KLS}}(x;\lambda)
=
{}_2F_0\!\left(
\begin{matrix}
-m,-x\\
-
\end{matrix};
-\lambda^{-1}
\right).
\]
Using the definition of the hypergeometric series
\citep[p.~5, (1.4.1)]{koekoek_lesky_swarttouw_2010}, which terminates
here because of the numerator parameter \(-m\), together with the
shifted-factorial identities on p.~4, gives
\[
C_m^{\mathrm{KLS}}(x;\lambda)
=
\sum_{j=0}^{m}
(-1)^j\binom{m}{j}
\frac{x(x-1)\cdots(x-j+1)}{\lambda^j}.
\]
where the product for \(j=0\) is understood as \(1\). Consequently,
\[
C_m(\lambda,x)
=
(-1)^m C_m^{\mathrm{KLS}}(x;\lambda).
\]
Thus, the formula in the preceding corollary is the classical
finite-sum expression written in the normalization used throughout
this paper.

\begin{proof}[Proof of Lemma \ref{lem:charlier-inversion}]
	Note the trivial identity
\[
\bigl(e^{wR} \left( w+1\right)^{-n}\bigr)\bigl(e^{-wR} \left( w+1\right)^n\bigr)=1.
\]
Recognizing both factors on the left-hand side as generating functions of Charlier polynomials~\eqref{eq:charlier-generating-function} with appropriate parameters, our identity takes the form
\[
\sum _{j\geqslant 0}\frac{(-R)^jC_{j}\left( -R ,-n\right)}{j!} w^j\sum _{m\geqslant 0}\frac{R^mC_{m}\left( R ,n\right)}{m!}w^m =1.
\]
Comparing the coefficients of \(w^s\) in the Taylor expansions on both
sides completes the proof of the lemma.
\end{proof}

For \(x\in\mathbb{R}\) and \(j\in\mathbb{N}_0\), define the falling factorial by
\[
x_{(0)}:=1,
\qquad
x_{(j)}:=x(x-1)\cdots(x-j+1),
\quad j\geqslant 1.
\]

\begin{lem}
\label{lem:charlier-tau-duality}
For all real values of \(t, y\) and every nonnegative integer \(N\), the following two identities hold, with \(t^NC_N(t,y)\) at \(t=0\) understood by polynomial continuation:
\begin{align}
\label{eq:charlier-tau-expansion-in-y}
t^{N} C_{N}(t, y)
&= \sum_{j=0}^{N} \binom{N}{j} \tau_{N-j}(y)\,(y-t)^{j}, \\[1mm]
\label{eq:charlier-tau-expansion-in-t}
t^{N} C_{N}(t, y)
&= \sum_{j=0}^{N} \binom{N}{j} \tau_{N-j}(t)\,(y-t)_{(j)}.
\end{align}
\end{lem}

\begin{proof}
For \(t\ne0\), the generating function \eqref{eq:charlier-generating-function} of the Charlier polynomials gives
\[
\sum_{m\geqslant0}\frac{t^{m}C_m(t,y)}{m!}w^{m} = e^{-wt}(1+w)^{y}.
\]
To derive~\eqref{eq:charlier-tau-expansion-in-y}, rewrite this as
\[
e^{-wt}(1+w)^{y} = \bigl(e^{-w}(1+w)\bigr)^{y} e^{w(y-t)}.
\]
The first factor \(\bigl(e^{-w}(1+w)\bigr)^{y}\) on the right according to \eqref{eq:tau-generating-function-appendix} is the exponential generating function of \(\tau_m(y)\),
and the second factor \(e^{w(y-t)}\)  is an exponential generating function of \((y-t)^{j}\). Hence,
\[
\sum_{m\geqslant0}\frac{t^{m}C_m(t,y)}{m!}w^{m}
= \left(\sum_{m=0}^{\infty}\frac{\tau_m(y)}{m!}w^{m}\right)
  \left(\sum_{j=0}^{\infty}\frac{(y-t)^{j}}{j!}w^{j}\right).
\]
Comparing coefficients of \(w^{N}\) on both sides gives~\eqref{eq:charlier-tau-expansion-in-y}.

Similarly, for~\eqref{eq:charlier-tau-expansion-in-t}, note the equivalent decomposition
\[
e^{-wt}(1+w)^{y} = \bigl(e^{-w}(1+w)\bigr)^{t}(1+w)^{y-t}.
\]
Here the first factor is the generating function of \(\tau_m(t)\), 
and the second generates the falling factorials \((y-t)_{(j)}\). 
Extracting the coefficient of \(w^{N}\) yields~\eqref{eq:charlier-tau-expansion-in-t}. Since both right-hand sides are polynomials in \(t\), the identities extend to \(t=0\), completing the proof.
\end{proof}
\begin{proof}[Proof of Lemma \ref{lem:charlier-pointwise-bound}]
Applying Proposition \ref{prop:tau-bound} to estimate the
\(\tau_{N-j}(y)\) terms on the right side of
identity~\eqref{eq:charlier-tau-expansion-in-y} in
Lemma~\ref{lem:charlier-tau-duality}, we get
\[
\begin{split}
	\bigl|t^NC_N(t,y)\bigr|
&\leqslant\sum ^{N}_{j=0}\binom{N}{j}\bigl|\tau_{N-j}(y)\bigr|| y-t|^{j}
\\
&\leqslant\sum ^{N}_{j=0}\binom{N}{j}(N-j)!|y|^{\frac{N-j}{2}}| y-t|^{j}
\\
&=N!|y|^{\frac{N}{2}}\sum ^{N}_{j=0}\frac{1}{j!} \left|\frac{y-t}{\sqrt{|y|}}\right|^{j}
\\
&\leqslant N!|y|^{\frac N2}e^{\frac{|y-t|}{\sqrt{|y|}}}.
\end{split}
\]
This completes the proof.
\end{proof}
\begin{proof}[Proof of Lemma \ref{lem:charlier-integral-bound}] By Lemma \ref{lem:charlier-pointwise-bound},
	\begin{equation}
	\label{eq:charlier-integral-reduction}
	\int_0^\infty t^{n} e^{-t}\bigl|t^NC_N(t,n)\bigr|\,dt\leqslant N!n^{N/2}\int_0^\infty t^{n} e^{-t}e^{\frac{|n-t|}{\sqrt{n}}}\,dt.		
	\end{equation}
With $u=(n-t)/\sqrt n$ we use $e^{|u|}< e^{u}+e^{-u}$ to bound the integral on the right:
	\[
	\begin{split}
		\int_0^\infty t^{n} e^{-t}e^{\frac{|n-t|}{\sqrt{n}}}\,dt
&\leqslant \int_0^\infty t^{n}e^{-t} \left(e^{\frac{n-t}{\sqrt{n}}}+e^{-\frac{n-t}{\sqrt{n}}}\right)\,dt
\\
&=n!\left(\frac{e^{\sqrt{n}}}{\left(1+\frac{1}{\sqrt{n}}\right)^{n+1}}+\frac{e^{-\sqrt{n}}}{\left(1-\frac{1}{\sqrt{n}}\right)^{n+1}}\right)
\\
&\leqslant 12n!.
	\end{split}
	\]
Here we used the inequality
\[
\frac{e^{1/x}}{(1+x)^{1/x^2}}<3
\]
that is valid for all \(x\geqslant-1/\sqrt{2}\) with \(x\ne0\), applying it with \(x=1/\sqrt{n}\) and \(x=-1/\sqrt{n}\) where \(n \geqslant 2\).
	Using this estimate to evaluate the integral on the right side of the inequality (\ref{eq:charlier-integral-reduction}) we complete the proof of the lemma.

\end{proof}

\begin{proof}[Proof of Lemma~\ref{lem:charlier-tau-shift}]
By identity \eqref{eq:charlier-tau-expansion-in-t} in Lemma~\ref{lem:charlier-tau-duality}, the term
\(\left( -R \right)^s C_s(-R, -n)\) can be expanded in terms of \(\tau_m\).  
Substituting this into the left-hand side of the lemma's identity gives
\begin{equation*}\begin{split}
\sum_{s=0}^j \frac{(-R)^s C_s(-R, -n)}{s!} \, \Delta_x^s x^j 
&= \sum_{s=0}^j \frac{\Delta_x^s x^j}{s!}  \sum_{m=0}^s \frac{s!}{(s-m)!m!}{\tau_{m}(-R)(R-n)_{(s-m)}}  \\
&= \sum_{m=0}^j \frac{\tau_m(-R)}{m!} \sum_{s=m}^\infty \frac{\Delta_x^{s } x^j}{(s-m)!} (R - n)_{(s-m)} .
\end{split}
\end{equation*}
The upper range of the inner sum is taken to be \(s=\infty\) instead of \(s=j\) since \(\Delta^sx^j=0\) when \(s>j\).
Taking here \(x=0\) and making a change of variables \(s\to s+m\) in the inner sum on the right, we get
\begin{equation} \label{eq:charlier-tau-newton-intermediate}
\begin{split}
\sum_{s=0}^j \frac{(-R)^s C_s(-R, -n)}{s!} \, \Delta_x^s x^j \big|_{x = 0}
&= \sum_{m=0}^j \frac{\tau_m(-R)}{m!} \sum_{s=0}^\infty \frac{\Delta_x^{s + m} x^j\big|_{x = 0}}{s!} (R - n)_{(s)} .
\end{split}
\end{equation}
The inner sum can be evaluated using Newton's interpolation formula (see Problem~6, p.~221 of \cite{comtet_1974}):
\[
P(y) = \sum_{s=0}^\infty \frac{\Delta_x^s P(x)\bigl|_{x=0}}{s!} \, y_{(s)},
\]
which holds for any polynomial \(P(x)\). Applying this with \(P(x) = \Delta^m x^j\) and \(y = R-n\) yields
\[
\Delta_x^m x^j \big|_{x = R - n} 
= \sum_{s=0}^\infty \frac{\Delta_x^{s + m} x^j \big|_{x = 0}}{s!} (R - n)_{(s)}.
\]
The right-hand side of this identity is exactly the inner sum on the right-hand side of \eqref{eq:charlier-tau-newton-intermediate}; therefore, replacing it by \(\Delta_x^m x^j \big|_{x = R-n}\) completes the proof of the lemma.
\end{proof}
\begin{proof}[Proof of Lemma \ref{lem:mahler-tau-stirling}]
	Let us recall (see \cite{comtet_1974}) the generating function for the Stirling numbers of the second kind:
 \[
 \sum ^{\infty }_{j=s}\Stirling{j}{s}\frac{z^{j}}{j!}=\frac{\left( e^{z}-1\right) ^{s}}{s!},
 \]
 Multiplying both sides by \(\tau_s(x)\) and summing over \(s\), we get
 \[
 \sum_{s=0}^\infty\tau_s(x)\sum ^{\infty }_{j=s}\Stirling{j}{s}\frac{z^{j}}{j!}= \sum_{s=0}^\infty\frac{\tau_s(x)}{s!}\left( e^{z}-1\right) ^{s}.
 \]
 Changing the order of summation on the left and expressing the sum on the right using the generating function \eqref{eq:tau-generating-function-appendix} of \(\tau_s\) we obtain the identity
 \[
  \sum_{j=0}^\infty\frac{z^{j}}{j!}\sum ^{j }_{s=0}\tau_s(x)\Stirling{j}{s}= \left((1+(e^z-1))e^{-(e^z-1)}\right)^{x}=e^{-x(e^z-1-z)}=\sum_{j=0}^\infty\frac{\rho_j(x)}{j!}z^j
 \]
 Hence, comparing the coefficients of \(z^j\) on both sides gives the identity of the lemma.
\end{proof}

\section{Examples of functions for which the asymptotic expansions converge}
\setcounter{equation}{0}
\label{app:convergent-expansion-examples}
In this appendix we present two nontrivial classes: a class of
sequences satisfying identity~\eqref{eq:ramanujan-inverted-series} and a class of
functions satisfying identity~\eqref{eq:ramanujan-convergent-identity}.

\begin{prop}

\label{prop:laplace-sequence-inversion}
Let \(F\) be a bounded nondecreasing function and define
\[
A_n=\int_{0}^{\infty}e^{-nt}\,dF(t),
\qquad n=0,1,2,\ldots.
\]
Then identity~\eqref{eq:ramanujan-inverted-series} holds for every nonnegative integer
\(n\).
\end{prop} 
\begin{proof}
Since \(F\) is bounded and nondecreasing, \(dF\) is a finite
positive measure. For \(x\geqslant0\), all terms obtained after inserting the
integral representation of \(A_m\) into the defining series for
\(f(x)\) are nonnegative. Hence Tonelli's theorem justifies the
exchange of summation and integration in the following calculation:
\begin{equation}
\label{eq:laplace-poisson-transform}
	f(x) = e^{-x}\sum_{m=0}^{\infty}\frac{A_{m}}{m!}\,x^{m}
     = e^{-x}\sum_{m=0}^{\infty}\frac{1}{m!}\int_{0}^{\infty} (x e^{-t})^{m}\, dF(t)= \int_{0}^{\infty} e^{-x(1-e^{-t})}\, dF(t).
\end{equation}
 On the other hand,  
\[
\Delta^{m}A_{n} = (-1)^{m}\int_{0}^{\infty} (1-e^{-t})^{m} e^{-nt}\, dF(t).
\]
Furthermore, Proposition~\ref{prop:tau-bound}, together with the endpoint identities \(\tau_0=1\) and \(\tau_m(0)=0\) for \(m\geqslant1\), implies that
\[
(-1)^m\tau_m(-n)\geqslant 0,
\qquad m,n\geqslant 0.
\]
Since \(1-e^{-t}\geqslant0\), all terms in the series occurring below are
nonnegative. Tonelli's theorem therefore justifies the second exchange
of summation and integration.
Recalling the generating function (\ref{eq:tau-generating-function-intro}) of the polynomials \(\tau_m\), we obtain
\[
\begin{aligned}
\sum_{m=0}^{\infty}\frac{\Delta^{m}A_{n}}{m!}\,\tau_{m}(-n) 
   &= \int_{0}^{\infty}\left(\sum_{m=0}^{\infty}\frac{\tau_{m}(-n)(-1)^{m}}{m!}(1-e^{-t})^{m}\right) e^{-nt}\, dF(t) \\
   &= \int_{0}^{\infty} e^{-n(1-e^{-t})}\, dF(t) \end{aligned}
\]
Since the last expression coincides with the expression \eqref{eq:laplace-poisson-transform} for \(f(n)\),  identity~(\ref{eq:ramanujan-inverted-series}) holds for the sequence \((A_n)\).

\end{proof}

\begin{prop}
\label{prop:moment-generating-ramanujan-identity}
Consider
\[
\varphi(x) = \int_{0}^{\infty} e^{xt}\, dF(t),
\]
where \(F(x)\) is a non-decreasing bounded function and \(R>0\) is such that
\[
\int_{0}^{\infty} e^{R(e^t-1)}\, dF(t) < \infty.
\]
For this choice of \(\varphi\) and \(R\), the identity (\ref{eq:ramanujan-convergent-identity}) holds.
\end{prop}
\begin{proof}
Let \(dF\) denote the finite positive measure induced by \(F\). The
moment assumption also justifies differentiation under the integral
sign. Indeed, for each \(m\geqslant0\) and \(|x-R|\leqslant R/2\),
\[
t^m e^{xt}\leqslant C_m e^{R(e^t-1)},
\qquad t\geqslant0,
\]
because the ratio of the left-hand side to the exponential on the
right is bounded on \([0,\infty)\). Dominated differentiation therefore
gives
\[
\varphi^{(m)}(R)=\int_0^\infty t^m e^{Rt}\,dF(t).
\]
All summands in the poissonization series are nonnegative, so Tonelli's
theorem justifies the exchange below:
\[
e^{-R}\sum_{m=0}^{\infty}\frac{\varphi(m)}{m!}R^{m}
   = e^{-R}\int_{0}^{\infty} \left(\sum_{m=0}^{\infty}\frac{(Re^t)^m}{m!}\right) dF(t)
   = \int_{0}^{\infty} e^{R(e^t-1)}\, dF(t).
\]

On the other hand, \(\rho_m(-R)\geqslant0\) for \(R>0\), by the
coefficient nonnegativity established above. Thus Tonelli's theorem
also justifies exchanging the Mahler-polynomial series and the integral.
Recalling the generating function~\eqref{eq:mahler-generating-function}, we have
\[
\sum_{m=0}^{\infty} \frac{\varphi^{(m)}(R)}{m!}\rho_m(-R)
   = \int_{0}^{\infty}\left(\sum_{m=0}^{\infty}\frac{\rho_m(-R)}{m!}t^m\right) e^{Rt}\, dF(t)
   = \int_{0}^{\infty} e^{R(e^t-1-t)} e^{Rt}\, dF(t).
\]
Comparing the right-hand sides of the last two expressions, we conclude that
identity (\ref{eq:ramanujan-convergent-identity}) is valid for \(\varphi\) and \(R\).
\end{proof}

\end{appendices}

\end{document}